\documentclass[12pt]{article} 
\usepackage[T2A]{fontenc}
\usepackage[cp1251]{inputenc}
\usepackage[english]{babel}
\usepackage{amsmath}
\usepackage{amsfonts}
\usepackage{amssymb}
\usepackage{hyperref}

\newcommand{\de}{\delta}
\newcommand{\eps}{\varepsilon}
\newcommand{\qed}{$\blacksquare$}
\newcommand{\ld}{\ldots}
\newcommand{\ges}{\geqslant}
\newcommand{\les}{\leqslant}
\newcommand{\sbs}{\subset}

\newcommand{\sign}{\mathrm{sign}}

\newcommand{\CC}{{\mathbb C}}
\newcommand{\NN}{{\mathbb N}}
\newcommand{\RR}{{\mathbb R}}
\newcommand{\ZZ}{{\mathbb Z}}

\newtheorem{Lem}{Lemma}
\newtheorem{Predl}{Proposition}

\newtheorem{Teor}{Theorem}
\newtheorem{TeorB}{Theorem} 
\newtheorem{Opr}{Definition}
\newtheorem{Zam}{Remark}
\newtheorem{Prim}{Example}

\author{Oleg~E.~Galkin$^1$, Svetlana~Yu.~Galkina$^2$}
\title{Functions consistent with real numbers,\\
and global extrema of functions\\ 
in exponential Takagi class}

\begin{document}

\maketitle

\footnotetext[1]{$^{,2}$National Research University Higher School of Economics, Nizhny Novgorod Branch.
Russia, Nizhny Novgorod, ul. Bolshaya Pecherskaya, 25/12.
E-mail$^{1}$: olegegalkin@ya.ru;
e-mail$^{2}$: svetlana.u.galkina@mail.ru}

\tableofcontents

\begin{abstract} 
The functions in exponential Takagi class are similar in construction to the continuous, 
nowhere differentiable Takagi function described in 1901.
They have one real parameter $v\in (-1;1)$ and for any real point $x$ are defined by the series 
$T_v(x) = \sum_{n=0}^\infty v^n T_0(2^nx)$, 
where $T_0(x)$ is the distance between $x$ and the nearest integer point.
If $v=1/2$ then $T_v$ coincides with Takagi function.
In this paper, for different values of the parameter $v$, we study the global extremums of the functions $T_v$, 
as well as the sets of extreme points.
All functions $T_v$ have a period $1$, so they are investigated only on the segment $[0;1]$.
This study is based on the properties of so called consistent and anti-consistent functions (polynomials or series), 
which the first half of the work is devoted to.
\end{abstract}




\bigskip
{\bf Preface.} 
The Russian version of the proposed paper was written in 2016 and edited in 2018, but was not published.
This English version is a draft of a future article.
It is just a translation of the 2018 text and does not yet take into account those results, 
which were independently obtained by other authors on similar topics in 2016 or later 
(see, for example, \cite{HanSchied2020}, \cite{MishuraSchied2019}, \cite{Han2019} and references in them). 
We plan to analyze mentioned results in the next versions of our paper.

\section{Introduction} 
\label{SectVved}
The work is devoted to the search for global extremes, as well as the points where they are reached, 
for functions from one single-parameter family containing the Takagi function. 
We call this family exponential Takagi class.
Our research is based on the properties of the so-called consistent and anti-consistent unitary polynomials and series, 
which are covered in the first half of the work.

Teiji Takagi defined his function $T(x)$ in~1901 in~\cite{Takagi}, 
where he also proved that $T(x)$ on $\mathbb R$ is everywhere continuous and nowhere differentiable. 
This function can be set by the series
$$
T(x) = \sum\limits_{n=0}^\infty \frac1{2^n} T_0(2^nx), \quad x\in {\mathbb R},
$$
where $T_0(x) = |x-[x+1/2]| = |\{x+1/2\}-1/2| = \rho(x, {\mathbb Z})$
is the distance between the point $x$ and the nearest integer point, 
$[y]$ is the integer part (respectively $\{y\}$ is the fractional part) of the number $y\in\RR$.
Obviously, the equalities $T_0(x)=x$ for $x\in[0;1/2]$ and $T_0(x)=1-x$ for $x\in[1/2;1]$ are correct.

Hata and Yamaguchi~\cite[Sec.~2]{HataYamaguti84} replaced the sequence of coefficients $\{1/2^n\}$ 
in the definition of the Takagi function with an arbitrary sequence of constants $\{c_n\}$ 
and got a new family of functions, calling it {\it Takagi class}.
The object of our study is the real functions $T_v(x)$ belonging to a narrower family.
They have on the parameter $v$ and are defined by the equality
\begin{equation} \label{eqDefTv}
T_v(x) = \sum\limits_{n=0}^\infty v^n T_0(2^nx), \quad x\in {\mathbb R}.
\end{equation}
Since the coefficients $c_n = v^n$ in this formula depend on $n$ according to the exponential law, 
then we will call the set of functions of the form~\eqref{eqDefTv}, where $v\in(-1;1)$, {\it exponential Takagi class}.

Note that for $v=0$, the function $T_v(x)$ matches $T_0(x)$, and for $v=1/2$ it matches with the Takagi function $T(x)$.

Obviously, functions $T_v$ have a period of 1 on $\mathbb R$ for all $v\in(-1;1)$, 
so we study their properties only on the segment $[0;1]$.

For each $v\in(-1;1)$, we denote by $M_v$ the global maximum of the $T_v$ function on the segment $[0;1]$, 
and denote by $E_v$ the set of points from $[0;1]$, where it is reached.

Throughout this work, the maximum, minimum, or extremum will be understood in a global sense and taken along the segment $[0;1]$.

\emph{Binary rational numbers} are numbers of the form $x=p/2^k$, where $p\in\mathbb Z$, $k\in {\mathbb N}\cup\{0\}$.

In this article, we find global extremes of the functions $T_v$ on the segment $[0;1]$ and the set of points of global extremes for all values of the parameter $v\in[0;1]$.
For this purpose, in the first half of the work, we define and study the properties of so-called consistent and anti-consistent unitary polynomials and series.

We can say that the problem of finding extremes for $T_v$ was set by J.~ Tabor and J.~ Tabor (see~\cite[Problem 1.2, p.~731]{Tabor2}.
In order to accurately estimate continuous semi-convex functions, they introduced the functions $\omega_p$, which can be set as $\omega_p(x) = 2\cdot T_{1/2^p}(x)$, and obtained the formula for global maxima of the functions $\omega_{p_n}$ on $[0;1]$ 
for one specific sequence $\{p_n\}$. 
Given our notations, their result~\cite[Theorem~3.1]{Tabor2} can be represented as follows:
\begin{TeorB}[J.~Tabor and J.~Tabor, 2009]
\label{TeorTabor}
Let $v_n$ be the unique positive solution of the equation $2v + 4v^2 +\ldots + 2^n v^n = 1$ 
for each $n\in\NN$, and $c_v = 1 / (1 - (4v-1)^{\log_{2v}v})$ for $v\in(1/4;1/2)$. 
Then 

a). $v_1=1/2$, the sequence $\{v_n\}$ decreases and converges to $1/4$.

b). $\max\limits_{x\in[0;1]} T_{1/4}(x) = 1/2 = \lim\limits_{v\to1 / 4+0} C_v$;\quad
$\max\limits_{x\in[0;1]} T_{1/2}(x) = 2/3 = \lim\limits_{v\to1/2-0} C_v$.

c). $\max\limits_{x\in[0;1]} T_{v_n}(x) = C_{v_n}$ for $n\in\NN$, $n\geqslant2$.
\end{TeorB}

Here are some more results about the extremes of the $T_v$ functions.

1) In the case of $v=1/2$ Kahane (see~\cite[Lieu~1]{Kahane}) found points of local and global extremes for $T_v (x)=T(x)$. In particular, he proved the following statement:
\begin{TeorB}[Kahane, 1959] 
\label{TeorKahane}
The set of points where the Takagi function $T(x)$ reaches its global maximum 
is the set of points that have a binary expansion $x=\ldots,x_1 x_2\ldots x_n\ldots$ 
and satisfy the condition $x_{2k+1}+x_{2k+2}=1$ for $k=0,1,\ldots$.
\end{TeorB}

2) Further results on local extremes and sets of the Takagi function level can be found in the reviews~\cite{AllaartKawamuraSurv} and~\cite{Lagarias}. See~also~\cite{Allaart}.

3) For $v=1/\sqrt{2}$, it follows from the results of~\cite[Lemma~5]{Galkina} that the global maximum of the $T_v$ function by $[0;1]$ is $(2+\sqrt{2})/3$ and is reached only at points $1/3$ and $2/3$.

4) In~\cite[Theorem 4]{Galkin2015} it is proved that for $v\in[-1/2;1/4]$, the point $1/2$ is the point of the global maximum of the function $T_v$ on $[0;1]$,
and for $v\in(-1;-1 / 2)\cup(1/4;1)$ is not.

5) In the case of $v=-1/2$, it follows from Allaart's remark~\cite[Remark 5.6, p.28]{Allaart} 
that the set of minimum points of the function $T_v$ is a Cantor type set obtained by removing the <<middle half>>, 
and therefore is uncountable.

The Takagi function and its generalizations are used in various fields of mathematics: mathematical analysis, probability theory, number theory, and others. A large number of publications are dedicated to these functions, 
and this number continues to increase. Lots of interesting results and links 
available in the reviews~\cite{AllaartKawamuraSurv} and~\cite{Lagarias}.

One of the applications of the Takagi function in number theory is its use in the Trollope formula (see~\cite[Theoreme~1]{Trollope}), which can be written as:
$$
\frac{1}{N}\sum_{n=1}^{N-1}s(n) - \frac{\log_2{N}}{2N} = 
\frac{1-\{\log_2{N}\}}{2} - \frac{T( 2^{ \{\log_2 N\}-1 } )}{2^{\{\log_2 N\} }},
$$
where $N\geqslant2$, and $s(n)$ is the sum of the binary digits of the number $n$.
Kruppel in~\cite[eq.~2.13]{Kruppel2008} found the global extremes of the right side of this formula.
This allowed him to get an accurate estimate for the left side
(previously such estimate was obtained in other ways, see the links in~\cite{Kruppel2008}).

Let us briefly describe the structure and main results of the work.
The work consists of six sections.
{\it Section~\ref{SectVved}} is an introduction.
{\it Section~\ref{SectIzvSvoistva}} contains reminder about some well-known properties of functions
in exponential Takagi class, given in~\cite{Galkin2015}.
{\it In~section~\ref{SectSoglFunc}} polynomials and series, consistent and anti-consistent with real numbers are introduced, 
and their properties are studied.
{\it In~section~\ref{SectSvazSoglSEkst}} usage of consistent polynomials and series to finding the global maxima of $T_v$ functions, 
and usage of anti-consistent polynomials and series to finding their global minima are described.
{\it In~section~\ref{SectGlobMin}} the results of the previous sections are used to calculate the global minima of the $T_v$ functions and the sets where they are reached.
Similarly, {\it in the last section~\ref{SectGlobMaks}} the global maxima of the $T_v$ functions, as well as the sets where they are reached, are studied and calculated, when possible.
At the end of this section, there are several examples of searching for global maxima.

\section{Some important properties of functions in exponential Takagi class}
\label{SectIzvSvoistva}

In this section, we present some properties of the $T_v$ functions that are necessary for further research 
(see \cite{AllaartKawamuraSurv}, \cite{Lagarias}, \cite{Galkin2015} and links in them).

\medskip
{\bf Domain of definition, evenness and continuity}

For any $v\in(-1;1)$, the function $T_v$ on $\RR$ is obviously even and satisfies the identity $T_v(1-x) = T_v(x)$.
In addition, in~\cite[Theorem~1]{Galkin2015} it is shown that

1) if $|v| < 1$, then the series~\eqref{eqDefTv} that sets the function $T_v(x)$, 
converges uniformly on $x$, its sum $T_v(x)$ is continuous, and $|T_v(x)| \leqslant 1/(2-2|v|)$ for all $x\in\RR$;

2) if $|v| \geqslant 1$, then the series~\eqref{eqDefTv} converges if and only if $x$ is a binary rational point.
In this case, the function $T_v (x)$ is discontinuous everywhere on the set of binary rational points.

\medskip
{\bf\boldmath Functional equation for $T_v(x)$}

From the formula~\eqref{eqDefTv}, which defines the function $T_v(x)$, it immediately follows that for all $v\in(-1;1)$, $N\in\NN$ and $x\in\RR$, the equality is met
\begin{equation} \label{eqFunkUr}
T_v(x) = \sum\limits_{n=0}^{N-1} v^n T_0(2^nx) + v^{N}T_v(2^{N}x).
\end{equation}

\medskip
{\bf\boldmath Explicit formula for $T_v$ when $v=1/4$}

Hata and Yamaguchi in~\cite[p.~335]{HataYamaguti83} showed that in the case of $v=1/4$, for all $x\in[0;1]$ the equality holds
$
T_v(x) = T_{1/4}(x) = 2(x-x^2).
$

\medskip
{\bf Differentiability}

1) For $v=1/2$, the function $T_v$ is nowhere differentiable on $\RR$ (see~\cite{Takagi}).

2) For $v=1/4$, this function is differentiable at all points of $x\in\RR$, except for integers (it follows from the explicit formula of Hata -- Yamaguchi).

3) For $v=0$, the function $T_v$ is differentiable at all points from $\RR$, except for half-integers.

4) For $1/2 \le |v| < 1$, the function $T_v$ is not differentiable at any point from $\RR$ (follows from Kono's theorem~\cite[Theorem 2]{Kono} on functions in Takagi class).

5) For all $|v| < 1/2$ except $v = 1/4$, the function $T_v$ is not differentiable at binary rational points and is differentiable at the other points of $\RR$ (see \cite[Theorem 2]{Kono}, \cite[Proposition 1.1.2]{Spurrier}, \cite[Theorem 3]{Galkin2015}).

\section{Unitary polynomials and series that are consistent and anti-consistent with numbers}
\label{SectSoglFunc}

The importance of studying consistent and anti-consistent polynomials and series that this section is devoted to, 
will be seen from the theorem~\ref{TeorPoiskGlobExtr}.
\subsection{Definitions of consistent and anti-consistent polynomials and series, and related notions}

\begin{Opr}
\label{OprUnit}
\emph{Polynomial} $c_0 + c_1 x +\ld+ c_n x^n$ or \emph{power series} $c_0 + c_1 x + c_2 x^2 + \ld$ 
is called \emph{unitary} if the free term $c_0$ equals $1$, and other coefficients equal either $-1$ or $1$.
Any polynomial can be considered as a power series with a finite number of non-zero coefficients.
\end{Opr}

\begin{Opr}
\label{OprLeksSravnShod}
1) \emph{Power series} $F(x) = c_0 + c_1 x + c_2 x^2 + \ld$  is called \emph{(lexicographically) less than power series} 
$G(x) = d_0 + d_1 x + d_2 x^2 + \ld$, if  there exists such $n\in\NN$ that $c_n < d_n$ and $c_k = d_k$ for all $k=0,1,\ld, n-1$.
This fact is denoted by $F \prec G$.

2) The convergence of power series is understood as their coefficient-wise convergence.
\end{Opr}

\begin{Opr}
\label{OprSoglFunc}
Let $w\in\RR$.

1) The unitary polynomial $P(x) = c_0 + c_1 x +\ld+ c_n x^n$ is called 
\emph{polynomial consistent with the point} $w$, 
and $w$ is called \emph{point consistent with the polynomial $P(x)$}, 
if $P(w)=0$ and for any $k=1,\ld,n$ the next inequalities are met:
\begin{equation} \label{eqSoglKoef}
c_k\cdot(c_0 + c_1 w +\ld+ c_{k-1} w^{k-1}) < 0.
\end{equation}
In this case, two unitary series $F_w^+(x)$ and $F_w^-(x)$, which have the form
\begin{equation} \label{eqF+}
\begin{array}{rl}
F_w^+(x) &= c_0+c_1 x+\ld+c_nx^n +c_0x^{n+1}+c_1x^{n+2}+\ld+c_nx^{2n+1} + \\
+ \ld &= P(x)(1+x^{n+1}+x^{2(n+1)}+\ld),
\end{array}
\end{equation}
and
\begin{equation} \label{eqF-}
\begin{array}{rl}
F_w^-(x) &= c_0+c_1 x+\ld+c_nx^n -c_0x^{n+1}-c_1x^{n+2}-\ld-c_nx^{2n+1} - \\
- \ld &= P(x)(1-x^{n+1}-x^{2(n+1)}-\ld),
\end{array}
\end{equation}
are called \emph{attached (to the polynomial $P$ and point $w$) series}.
\emph{Intermediate series for the polynomial $P$ and point $w$} are defined as any unitary series 
(including attached ones), that have the form
\begin{equation} 
\begin{array}{rl}
&(c_0+c_1 x+\ld+c_nx^n) \pm (c_0x^{n+1}+c_1x^{n+2}+\ld+c_nx^{2n+1}) \pm \ld = \\
&= P(x)(1 \pm x^{n+1} \pm x^{2(n+1)} \pm \ld),
\end{array}
\end{equation}
where any combinations of <<$+$>> and <<$-$>> characters are allowed.

2) The unitary series $F(x) = c_0 + c_1 x + c_2 x^2 +\ld$ is called \emph{series consistent with the point $w$}, 
and $w$ is called \emph{the point consistent with the series $F(x)$},
if the inequalities~\eqref{eqSoglKoef} are met for any $k\in\NN$.
In this case, \emph{an attached and an intermediate series for the point $w$}
are assumed to be equal to the consistent series.

3) \emph{A function consistent with the point $w$} is defined as either polynomial or series that is consistent with that point.
\end{Opr}

\begin{Zam}
\label{ZamSoglFun}
1) The consistent function $c_0 + c_1 x +\ld$ always has $c_0=1$ (due to the definition~\ref{OprUnit}) 
and $c_1=-1$ (due to the inequality~\eqref{eqSoglKoef} with $k=1$).

2) Two different points from $\RR$ can be consistent with the same function (see~theorem~\ref{TeorSoglW<1/2>1}).
However, different points from the segment $[1/2;1]$ are always consistent with different functions (see~theorem~\ref{TeorProizEdinSogl}).

3) For any $w\in\RR$, the next lexicographic inequalities hold for the consistent function $F_w$ and the attached series $F_w^\pm$: 
$F_w^-\preceq F_w \preceq F_w^+$ (by virtue of the definitions \ref{OprLeksSravnShod} and~\ref{OprSoglFunc}).
\end{Zam}

Definitions of anti-consistent functions are similar to the definitions of consistent functions:

\begin{Opr}
\label{OprASoglFunc}
Let $w\in\RR$.

1) The unitary polynomial $Q(x) = c_0 + c_1 x +\ld+ c_n x^n$ is called 
\emph{polynomial anti-consistent with the point $w$}, 
and $w$ is called \emph{point anti-consistent with the polynomial $Q(x)$}, 
if $Q(w)=0$ and for any $k=1,\ld,n$ the next inequalities are met:
\begin{equation} \label{eqASoglKoef}
c_k\cdot(c_0 + c_1 w +\ld+ c_{k-1} w^{k-1}) > 0.
\end{equation}

2) The unitary series $A(x) = c_0 + c_1 x + c_2 x^2 +\ld$ is called \emph{series anti-consistent with the point $w$}, 
and $w$ is called \emph{point anti-consistent with the series $A(x)$},
if the inequalities~\eqref{eqASoglKoef} are met for any $k\in\NN$.

3) \emph{A function anti-consistent with the point $w$} is defined as either polynomial or series that is anti-consistent with that point.
\end{Opr}

Now let us prove the existence and uniqueness of consistent and anti-consistent functions.
\begin{Predl}
For any real point, there exists a unique function (a polynomial or a series) that is consistent with it, 
as well as a unique function (a polynomial or a series) that is anti-consistent with this point.
\end{Predl}
{\bf Proof.} Let $w\in\RR$.

1) The coefficients of the function consistent with $w$ can be constructed by induction.
First, we put $c_0=1$. 
Further, if the coefficients $c_0, c_1, \ld, c_n$ have already been constructed, then we denote $S_n(x) = c_0 + c_1 x +\ld+ c_n x^n$ and consider two possible cases:

a) if $S_n (w)=0$, then $S_n (x)$ is the desired consistent polynomial;

b) if $S_n(w)\neq 0$, then put $c_{n+1} = - \sign(S_n(w))$.

As a result, either in a finite number of steps we get a polynomial that is consistent with the point $w$, 
or in an infinite number of steps we get a series that is consistent with the point $w$.

2) The construction of an anti-consistent function is similar, 
but in the case of $S_n (w)\neq 0$ you need to take $c_{n+1} = \sign(S_n(w))$.
\qed 

\medskip
In the following statement, we find out when rational numbers are (anti-)consistent with polynomials.
\begin{Predl}
\label{PredlSoglIrraz}
Rational numbers cannot be the roots of unitary polynomials 
(hence they cannot be consistent or anti-consistent with them), except in two cases: 
the number $1$ is consistent with the polynomial $1-x$, 
and the number $-1$ is anti-consistent with the polynomial $1+x$.
\end{Predl}
{\bf Proof.}
Let the rational number is represented by the irreducible fraction $p/q$,
and let $p/q$ be the root of the unitary polynomial $P(x) = c_0 + c_1 x +\ld+ c_n x^n$.
Then by~\cite [Theorem~6, p.~256 ]{Kostrikin} the numbers $p$ and $q$ are divisors for $c_0=1$ and $c_n=\pm1$, respectively. 
So $p/q=\pm1$. 
\qed 

\subsection{Formulae for anti-consistent polynomials and series for all real numbers}

In this section, we calculate anti-consistent functions for each point of the numeric axis.

\begin{Teor} 
\label{TeorVseAntisogl}

1). All points $w>-1$ are anti-consistent with the same series $1+x+x^2+\ldots$. 
If $|x|<1$, then its sum is $1/(1-x)$.

2). The point $w=-1$ is anti-consistent with the polynomial $1+x$. 

3). All points $w<-1$ are anti-consistent with the same series $1+x-x^2-x^3+x^4+x^5-\ldots$ having alternating pairs of signs. 
For $|x|<1$, its sum is $(1+x)/(1+x^2)$.
\end{Teor}
{\bf Proof.}
The anti-consistency of the point $w=-1$ with the polynomial $1+x$ is obvious.
In the other two variants, due to the definition~\ref{OprASoglFunc}, 
it is enough to check whether the next inequality is correct for any $k\in\NN$:
$$
D_k = c_k\cdot(c_0 + c_1 w +\ld+ c_{k-1} w^{k-1}) > 0.
$$

If $w>-1$, then, by condition of theorem, $c_k = 1$, $k=0,1,\ld$.
So $D_k = 1 + w +\ld+ w^{k-1}$.
For $w\ges0$, the inequality $D_k>0$ is obviously satisfied, and for $-1<w<0$ we have: $D_k = (1-w^k)/(1-w) > 0$.

If $w<-1$, then the following four cases are possible for the number $k$:

a). If $k=4i$, where $i\in\NN$, then $c_k=1$ and
$$
D_k = 1 + w - w^2 - w^3 + \ld-w^{4i-1} = (1+w)(1-w^2)(1+w^4+w^8+\ld+w^{4i-4}) > 0. 
$$

b). If $k=4i+1$, where $i\in\NN$, then $c_k=1$ and
$$
D_k = 1 + w - w^2 - w^3 + \ld - w^{4i-1} + w^{4i} = D_{4i} + w^{4i} > 0. 
$$

c). If $k=4i+2$, where $i\in\NN$, then $c_k=-1$ and
$$
D_k = - (1 + w - w^2 - w^3 + \ld - w^{4i-1} + w^{4i} + w^{4i+1}) =
$$
$$
= -(1+w)(1 - (1-w^2)(w^2+w^6+\ld+w^{4i-2})) > 0. 
$$

d). If $k=4i+3$, where $i\in\NN$, then $c_k=-1$ and
$$
D_k = - (1 + w - w^2 - w^3 + \ld - w^{4i-1} + w^{4i} + w^{4i+1} - w^{4i+2}) =
D_{4i+2} + w^{4i+2} > 0. 
$$

So, $D_k>0$ for all $k\in\NN$, that is what we needed to check. 
\qed

\subsection{Properties of consistent functions that are common for all real numbers}

\begin{Lem}[the criterion of consistency of the polynomial]
\label{LemKritSoglMnog}
Let polynomial $P(x)=c_0+c_1x+\ld+c_nx^n$ be unitary and $w\in\RR$.
Then $P(x)$ is consistent with the point $w\in\RR$ if and only if: 
$P(w)=0$ and for any $k=1,\ld, n$ the inequality is satisfied
$$ 
c_k\cdot(c_kw^k+c_{k+1}w^{k+1}+\ld+c_nw^n) > 0.
$$ 
\end{Lem}
{\bf Proof} 
follows from the definition~\ref{OprSoglFunc} of a consistent polynomial and from the equality
$$
\begin{array}{rl}
c_0 + c_1 w +\ld+ c_{k-1} w^{k-1} &= P(w) - (c_kw^k+c_{k+1}w^{k+1}+\ld+c_nw^n) = \\
&= -(c_kw^k+c_{k+1}w^{k+1}+\ld+c_nw^n).
\quad\quad\text{\qed}
\end{array}
$$

\begin{Predl}[on the analyticity of the sum of the consistent series] 
\label{PredlAnalF}
1) In the complex plane, the convergence region of any unitary series (including consistent, attached, and intermediate ones) is the unit circle $\{z\in\CC \bigm| |z|<1\}$. 

2) In the circle $\{z\in\CC \bigm| |z|<1\}$, the sums of unitary series are analytic functions. 
In particular, the sums of the attached series $F^+$ and $F^-$ for consistent polynomial $P$ of degree $N$,
have the form
\begin{equation} \label{eqSumF+-}
F^+(z) = P(z)/(1-z^{N+1}) \quad\text{and}\quad
F^-(z) = P(z)\cdot(1-2z^{N+1})/(1-z^{N+1}).
\end{equation}
\end{Predl} 
{\bf Proof.} 
By virtue of the Cauchy-Hadamard formula (see~\cite[ch.~II, \S3, p.~5, p.~63]{Privalov}), 
the radius of convergence of any unitary series $F(z) = c_0 + c_1 z + c_2 z^2 +\ld$ 
is $ 1/\overline {\lim}_{n\to\infty}\sqrt[n]{|c_n|} = 1$.
At any point $z\in\CC$ with $|z|=1$, such a series diverges, since the common member of the series does not tend to $0$.
So, the unit circle $\{z\in\CC \bigm| |z|<1\}$ is a convergence region for this series, 
and its sum $F(z)$ is analytic within this circle (see~\cite[ch.~II, \S4, p.~6, p.~77]{Privalov}).
The formulae~\eqref{eqSumF+-} follow from the formulae \eqref{eqF+} and \eqref{eqF-}.
\qed 

\begin{Lem}[about the number of identical coefficients]
\label{LemOdinZnak}
If $m\in\NN$, $w$ is positive and $1-w - \ld-w^m < 0$, then 
no more than $m$ of identical non-zero coefficients can go in a row
in the function (polynomial or series) which is consistent with the point $w$.
\end{Lem} 
{\bf Proof.} 
According to item~1) of remark~\ref{ZamSoglFun} we have: $c_0=1$, $c_1=-1$.
Therefore, it is sufficient to prove that after any block of $(m+1)$ coefficients 
having the form $c_n=-1$, $c_{n+1} = c_{n+2} =\ld= c_{n+m} = 1$, 
the following coefficient $c_{n+m+1}$ is equal to $-1$ (for a block with opposite signs, the proof is similar).
If there is such a block, the function that is consistent with the point $w$ has the form:
$$
F_w(x) = c_0 + c_1 x +\ld+ c_{n-1}x^{n-1} - x^n + x^{n+1} +\ld+ x^{n+m} + \ld.
$$
From the inequality~\eqref{eqSoglKoef} for $k=n$, we have:
$c_0 + c_1 w +\ld+ c_{n-1}w^{n-1} > 0$.
From here and from the inequality $1-w - \ld-w^m < 0$ we get:
$$
c_0 + c_1 w +\ld+ c_{n+m}w^{n+m} =
c_0 + c_1 w +\ld+ c_{n-1}w^{n-1} - w^n(1-w - \ld - w^m) > 0.
$$
Therefore, $c_{n+m+1} = -\sign(c_0 + c_1 w +\ld+ c_{n+m}w^{n+m}) = -1$.
\qed 

\begin{Predl}[on convergence of unitary functions]
\label{PredlShodUnitFunc}
If the sequence of unitary functions $F_m$, $m\in\NN$, converges coefficient-wise to the unitary function $F$ at $m\to\infty$, 
then for any $z\in\CC$, such that $|z|<1$, $F_m(z)$ converges to $F(z)$ at $m\to\infty$.
\end{Predl}
{\bf Proof.}
Let $|z|<1$, $F_m(x) = c_0(m)+c_1(m)x+c_2(m)x^2 +\ld$ when $m\in\NN$,
and $F(x) = d_0+d_1x+d_2x^2 +\ld$.
Then
$$
|F_m(z)-F(z)| \les \sum_{k=0}^\infty |c_k(m)-d_k|\cdot|z|^k
\les \sum_{k=0}^n |c_k(m)-d_k|\cdot|z|^k + 2\sum_{k=n+1}^\infty |z|^k.
$$
For any $\eps>0$, there exists some $n\in\NN$, such that the second term will be less than $\eps/2$, 
since the series $\sum_{k=0}^\infty |z|^k$ converges. 
Let's fix this $n$. 
Then for some $M\in\NN$, for all $m\ges M$ the first term will also be less than $\eps/2$, 
since $c_k(m)\to d_k$ for $m\to\infty$, $k=1,\ldots, n$.
So $|F_m(z)-F(z)|<\eps$ for all $m\ges M$.
\qed 

\subsection{\boldmath Consistent functions for all points satisfying the condition $w\les 1/2$ or $w \ges 1$}
The situation with consistent polynomials and series is more complicated than with anti-consistent ones.
In this section, we find consistent functions only for points that lie in the intervals 
$(- \infty; 1/2]$ or $[1;+\infty)$.
The points of the remaining interval $(1/2;1)$ are considered in the following subsections.

\begin{Teor}
\label{TeorSoglW<1/2>1}

1). All points $w>1$ are consistent with the same series $1-x+x^2-x^3+x^4 - \ld$, whose signs alternate. If $|x|<1$, its sum is $1/(1+x)$.

2). The point $w=1$ is consistent with the polynomial $1-x$.

3). All points $w\in[-1; 1/2]$ are consistent with the same series $1-x-x^2-x^3 - \ld$.
If $|x|<1$, its sum is $(1-2x)/(1-x)$.

4). Consistent functions for points $w < -1$ can be described as follows:

4a) For any $k\in {\mathbb N}$, the polynomial $P_{2k}(t) = 1-t-t^2 - \ld-t^{2k}$ has a single negative root $t=w_k$, and $w_k$ belongs to the interval $(-2;-1)$.
The sequence $\{w_k\}$ strictly increases and converges to $-1$, so that the next asymptotic equality is true:
\begin{equation} \label{eqwkAsimp}
w_k = -1-\ln3/(2k)+\underline O (1/k^2) \quad (k\to\infty).
\end{equation}

4b) If $k\in {\mathbb N}$, and $w\in(w_{k-1}, w_k)$ (where $w_0=-\infty$), then the point $w$ is consistent with the series
\begin{equation} \label{eqFwk}
F_k(x) = 1 - \sum_{n=1}^{2k} x^n + \sum_{n=0}^{\infty}(-1)^n(x^{2k+2n+1}+x^{2k+2n+2}).
\end{equation}
For $|x|<1$, its sum is $(1-2x)/(1-x) + 2x^{2k+1}/((1-x)(1+x^2))$.

4c) For any $k\in {\mathbb N}$, the point $w_k$ is consistent with the polynomial
\begin{equation} \label{eqPwk}
P_{2k}(x) = 1-x-x^2 - \ld-x^{2k}.
\end{equation}
\end{Teor}
{\bf Proof.}
Due to the definition~\ref{OprSoglFunc}, in each of the four cases, 
it is sufficient to check that the next relations are met for all $i\in\NN$:
\begin{equation} \label{eqDi<0}
D_i = c_i\cdot(c_0 + c_1 w +\ld+ c_{i-1} w^{i-1}) < 0.
\end{equation}

1). In the first case, $w>1$ and $c_i=(-1)^i$ for all $i\in\NN$, so 
$$
D_i = (-1)^i\cdot(1-w +\ld+ (-1)^{i-1}w^{i-1}) = (-w^i+(-1)^i)/(1+w) < 0.
$$

2). The consistency of the point $w=1$ with the polynomial $1-x$ is obvious.

3). Let $w\in[-1; 1/2]$, $c_0=1$, $c_i=-1$ for all $i\in\NN$. 
Then $1-w>0$ and $2w-1-w^i < 0$ 
(that is true due to the inequality $2w-1-w^i < -1-w^i$ when $-1 \les w < 0$,  
due to $2w-1-w^i\les 2w-1$ when $0\les w < 1/2$, and due to the equality $2w-1-w^i = -1/2^i$ when $w = 1/2$).
Hence $D_i = -(1-w - \ld - w^{i-1}) = (2w-1-w^i)/(1-w) < 0$.

4). Let $w < -1$ and $k\in\NN$. Next, let us prove items 4a)-4c).

4a) For any $t<0$, the equality $P_{2k}(t) = (1-2t+t^{2k+1})/(1-t)$ is met, 
so the polynomials $P_{2k}(t)$ and $Q_{2k+1}(t) = 1-2t+t^{2k+1}$ have the same negative roots.
The derivative $Q'_{2k+1}(t) = -2 +(2k+1) t^{2k}$ has roots $\pm t_k$, where $t_k = (2/(2k+1))^{1/(2k)}\in(0;1)$.
Therefore, $Q_{2k+1}$ grows when $t\les-t_k$ and decreases when $-t_k\les t\les 0$.
So $Q_{2k+1}(t) > Q_{2k+1}(-1) = 2 > 0$ when $-1\les t\les -t_k$, $Q_{2k+1}(t) > Q_{2k+1}(0) = 1 > 0$ when $-t_k\les t\les 0$.
Then $Q_{2k+1}(t)$ is positive when $-1\les t\les 0$ and increases when $t\les-1$.
Also, $Q_{2k+1}(-2) = 5-2^{2k+1} < 0$.
Therefore, $Q_{2k+1}$ (and hence $P_{2k}$) has exactly one negative root $w_k$, moreover $w_k\in(-2;-1)$.

Now let's show that $w_k<w_{k+1}$.
We have: $Q_{2k+1}(w_k) = Q_{2k+3}(w_{k+1}) = 0$.
Then
$$
Q_{2k+3}(w_k) = Q_{2k+1}(w_k) - w_k^{2k+1} + w_k^{2k+3} = 
w_k^{2k+3} - w_k^{2k+1} = w_k^{2k+1}(w_k^2-1) < 0. 
$$
Hence $w_k<w_{k+1}$, because $Q_{2k+3}(t)$ increases when $t\les-1$.

Let's move on to the proof of the formula~\eqref{eqwkAsimp}.
To do this, we find the asymptotic sequence $Q_{2k+1}(-1-\ln3/(2k)-a/k^2)$ for an arbitrary $a\in\RR$
up to members of the order of $1/k$:
$$
Q_{2k+1}(-1-\ln3/(2k)-a/k^2) = 1 - 2(-1-\ln3/(2k)-a/k^2) + 
$$
$$
+ (-1-\ln3/(2k)-a/k^2)^{2k+1} =
3 + \ln3/k - e^{(2k+1)\ln(1+\ln3/(2k)+a/k^2)} + \overline{o}(1/k) =
$$
$$
= (\ln3(3\ln3-2) - 24a)/(4k) + \overline{o}(1/k).
$$
From here we get, for sufficiently large values of $k$:

if $a=0$, then $Q_{2k+1}(-1-\ln3/(2k)) = \ln3(3\ln3-2)/(4k) + \overline{o}(1/k) > 0$, so $w_k < -1-\ln3/(2k)$;

if $a=1$, then $Q_{2k+1}(-1-\ln3/(2k)-1/k^2) = (\ln3(3\ln3-2)-24)/(4k) + \overline{o}(1/k) < 0$, so $w_k > - 1-\ln3/(2k) -1/k^2$.

Hence, the asymptotic equality~\eqref{eqwkAsimp} is satisfied.

4b) Here we show that all numbers $w\in(w_{k-1}, w_k)$ are consistent 
with the series $F_k(x) =\sum_{i=0}^{\infty}c_ix^i$ of the form~\eqref{eqFwk}.
In this case, $Q_{2k-1}(w)>0$, $Q_{2k+1}(w)<0$.
To check the relations~\eqref{eqDi<0} for $i\in\NN$, consider five sub-cases.

4b$_1$) Let $1\les i \les 2k$.
Then for odd $i = 2j-1 \les 2k-1$ we have:
$Q_{2j-1}(w) > Q_{2j-1}(w_{k-1}) \ges Q_{2j-1}(w_{j-1}) = 0$.
So
$$
D_i = D_{2j-1} = -(1-w - \ld-w^{2j-2}) = - Q_{2j-1}(w)/(1-w) < 0.
$$
Therefore, for even $i=2j\leqslant 2k$ we get:
$$
D_i = D_{2j} = -(1-w - \ld-w^{2j-1}) = D_{2j-1} + w^{2j-1} < 0.
$$
The remaining four sub-cases are devoted to the option $i > 2k$.

4b$_2$) Let $i = 2k+4j-3$, where $j\in\NN$.
Then $c_i=c_{2k+4j-3}=1$, and
$$
D_i = D_{2k+4j-3} = 1 - \sum_{n=1}^{2k} w^n + \sum_{n=0}^{2j-3}(-1)^n(w^{2k+2n+1}+w^{2k+2n+2}) =
$$
$$
= Q_{2k+1}(w)/(1-w) + w^{2k+1}(1+w)(1-w^{2(2j-2)})/(1+w^2) < 0.
$$

4b$_3$) Let $i = 2k+4j-2$, where $j\in\NN$.
Then $c_i=1$, and
$$
D_i = D_{2k+4j-2} = D_{2k+4j-3} + c_{2k+4j-3}w^{2k+4j-3} = D_{2k+4j-3} + w^{2k+4j-3} < 0.
$$

4b$_4$) Let $i = 2k+4j-1$, where $j\in\NN$.
Then $c_i=c_{2k+4j-1}=-1$, and
$$
D_i = D_{2k+4j-1} = -( 1 - \sum_{n=1}^{2k} w^n + \sum_{n=0}^{2j-2}(-1)^n(w^{2k+2n+1}+w^{2k+2n+2})) =
$$
$$
= -( 1 - \sum_{n=1}^{2k-2} w^n) - \sum_{n=-1}^{2j-2}(-1)^n(w^{2k+2n+1}+w^{2k+2n+2}) =
$$
$$
= Q_{2k-1}(w)/(w-1) + w^{2k-1}(1+w)(1-w^{2(2j)})/(1+w^2) < 0.
$$

4b$_5$) Let $i = 2k+4j$, where $j\in\NN$.
Then $c_i=-1$, and
$$
D_i = D_{2k+4j} = D_{2k+4j-1} - c_{2k+4j-1}w^{2k+4j-1} = D_{2k+4j-1} + w^{2k+4j-1} < 0.
$$

4c) For any $k\in{\mathbb N}$, the check of consistency of point $w=w_k$ with the polynomial~\eqref{eqPwk}
is the same as in item~4b$_1$) of this proof.
\qed 

\subsection{\boldmath Any point in the interval $[1/2;1)$ is the root of its consistent function}
It is not possible to get explicit formulae of consistent functions for all points $w\in[1/2;1)$.
For this reason, the current and following subsections of the this section 
are devoted only to analysis of some properties of these functions.

\begin{Teor}[on roots of the consistent series] 
\label{TeorFw(w)=0}
If the point $w$ belongs to half-open interval $[1/2;1)$ and is consistent with a series (not a polynomial), 
then $w$ is the root of this consistent series.
\end{Teor} 
{\bf Proof.} 
Let the series $F_w$ be consistent with the point $w$.
We must show that $F_w(w)=0$.
If $w=1/2$, then according to item~3) of the theorem~\ref{TeorSoglW<1/2>1} for $|x|<1$ we have: $F_{1/2}(x) = (1-2x)/(1-x)$.
So $F_{1/2}(1/2) = 0$.

If $w\in(1/2;1)$, then $F_{1/2}(w) = 1-w-w^2 - \ld = (1-2w)/(1-w) < 0$, so for some $m\in\NN$ it will be $1-w-w^2-\ld-w^m < 0$. 
Then, due to lemma~\ref{LemOdinZnak}, the series $F_w$ cannot have more than $m$ identical coefficients in a row.
Therefore, the sign of the coefficients of the series changes an infinite number of times. 
If we denote by $N_k$ the increasing numbers of terms, which changes its sign, 
then $c_{N_k} = - c_{N_k-1}$ for any $k\in\NN$.
So for any $\eps>0$, there exists a number $K\in\NN$, such that $w^{N_K}<\eps$.
It remains for us to show that for any $n>N_K$ the partial sum $S_n = c_0 + c_1 w +\ld+ c_n w^n$
satisfies the inequality $|S_n|<\eps$. 
From here the required equality will follow: $F_w(w) =\lim_{n\to\infty}S_n = 0$.

So let $n>N_K$. Then $N_p < n\les N_{p+1}$ for some $p\ges K$, so $w^{N_p}<\eps$.
Next, we have: $S_{N_p} = S_{N_p-1} + c_{N_p}w^{N_p}$.
First, consider the case of $c_{N_p}=1$. 
Then $c_{N_p+1} = -c_{N_p} < 0$ by definition $N_p$.
Therefore, due to the consistency of the series $F_w$, $S_{N_p} > 0$ and $S_{N_p-1}<0$.
So $0 < S_{N_p} < w^{N_p}$.
Similarly, we get that $ - w^{N_{p+1}} < S_{N_{p+1}} < 0$.
Further, since $N_p < n\les N_{p+1}$, then $c_n=-1$ and $S_n = S_{N_p} - w^{N_p+1} - \ld - w^n$.
Hence, on the one hand, $S_n < S_{N_p} < w^{N_p}$.
On the other hand, we have: $S_n \ges S_{N_p}-w^{N_p+1}-\ld-w^{N_{p+1}} = S_{N_{p+1}} > -w^{N_{p+1}}$.
So, $|S_n| < w^{N_p}<\eps$.

The similar reasoning is true for $c_{N_p}=-1$.
\qed 

\medskip
The following criterion follows from the proved theorem.
\begin{Lem}[the criterion of consistency of the series]
\label{LemKritSoglRad}
The unitary series $F(x)=c_0+c_1x+c_2x^2+\ld$ is consistent with the point $w\in[1/2;1)$ if and only if $F(w)=0$ and for any $k\in\NN$ the inequality is satisfied
$$
c_k\cdot(c_kw^k+c_{k+1}w^{k+1}+\ld) > 0.
$$
\end{Lem}
{\bf Proof} 
Given the theorem~\ref{TeorFw(w)=0}, the proof of this lemma is similar to the proof of the lemma~\ref{LemKritSoglMnog}.
\qed 

\subsection{Changing the consistent functions when extracting the square root or squaring}

\begin{Teor}[about extracting the square root]
\label{TeorSqrtC}
1) If $w\in[1/\sqrt2, 1)$, then the function that is consistent with the point $w$ has the form 
$F_w(x) = (1-x)F_{w^2}(x^2)$, where function $F_{w^2}$ is consistent with the point $w^2\in[1/2, 1)$.

2) If the point $w\in[1/2, 1)$ is consistent with the function $F_w$, then the point $\sqrt{w}\in[1/\sqrt2, 1)$ is consistent with the function $F_{\sqrt w}(x) = (1-x)F_w(x^2)$.

3) If $w\in[1/\sqrt[2^n]{2}, 1)$, then the function consistent with $w$ has the form
\begin{equation} \label{eqsqrt2nw}
F_w(x) = (1-x)(1-x^2)\cdots(1-x^{2^{n-1}}) F_d(x^{2^n}),
\end{equation}
where $F_d$ is the function consistent with the point $d = w^{2^n}\in[1/2, 1)$.
In this case, if $F_d$ is a polynomial, then $F_w$ is a polynomial too, and $\deg F_w = 2^n(\deg F_d + 1) - 1$.
\end{Teor}
{\bf Proof.} 
1) We assume that the function $F_{w^2}(x) = c_0+c_1x+\ld+c_kx^k+\ld$ which is consistent with the point $w^2$, 
is a series, since if $F_{w^2}(x)$ is a polynomial, then the proof is similar.
Let us show that the function $F(x) = (1-x)F_{w^2}(x^2)$ is consistent with the point $w$.
We have:
$$
F(x) = c_0-c_0x + c_1x^2-c_1x^3 + \ld + c_kx^{2k}-c_kx^{2k+1} + \ld = b_0+b_1x+b_2x^2+\ld,
$$
where $b_{2k} = c_k$, $b_{2k+1} = -c_k$ when $k=0,1,2,\ld$.
Therefore, the series $F(x)$ is unitary.
In addition, $F(w) = (1-w)F_{w^2}(w^2) = 0$ by the theorem~\ref{TeorFw(w)=0}.
Thus, by virtue of the lemma~\ref{LemKritSoglRad}, in order to prove the consistency of $F(x)$ with the point $w$, 
it remains to check the next inequalities for any $i\in\NN$:
\begin{equation} \label{eqbicdotgt0}
R_i = b_i\cdot(b_iw^i+b_{i+1}w^{i+1}+\ld) > 0.
\end{equation}

If $i=2k$, we have:
$$
R_i = b_{2k}\cdot(b_{2k}w^{2k}+b_{2k+1}w^{2k+1}+\ld) = 
$$
$$
= c_k\cdot(c_kw^{2k}-c_kw^{2k+1}+c_{k+1}w^{2k+2}-c_{k+1}w^{2k+3}+\ld) =
$$
$$
= (1-w)\cdot c_k\cdot(c_k(w^2)^k+c_{k+1}(w^2)^{k+1}+\ld).
$$
Since the series $F_{w^2}(x)$ is consistent with $w^2$, then 
$c_k\cdot(c_k(w^2)^k+c_{k+1}(w^2)^{k+1}+\ld) > 0$ by lemma~\ref{LemKritSoglRad}.
Hence, $R_i>0$.

If $i=2k+1$, we have:
$$
R_i = b_{2k+1}\cdot(b_{2k+1}w^{2k+1}+b_{2k+2}w^{2k+2}+\ld) = 
$$
$$
= - c_k\cdot(-c_kw^{2k+1}+c_{k+1}w^{2k+2}-c_{k+1}w^{2k+3}+c_{k+2}w^{2k+4}-c_{k+2}w^{2k+5}+\ld) =
$$
$$
= w^{2k+1}\cdot(1 - c_k\cdot(1-w)\cdot(c_{k+1}w+c_{k+2}w^3+\ld)) \ges 
$$
$$
\ges w^{2k+1}\cdot(1 - (1-w)\cdot(w+w^3+\ld)) = w^{2k+1}/(1+w) > 0.
$$
So, the inequality~\eqref{eqbicdotgt0}, and the item~1) of the theorem are proved.

2) The statement of item~2) of the theorem can be obtained from the statement of item~1) 
by replacing the number $w$ with $\sqrt w$.

3) The formula~\eqref{eqsqrt2nw} can be proved by $n$--multiple application of the statement of item~1).
This formula implies the equality
$$
\deg F_w = 1+2+\ld+2^{n-1} + \deg F_d\cdot 2^n = 2^n(\deg F_d + 1) - 1.
$$
The theorem~\ref{TeorSqrtC} is proved. \qed 

\medskip
In the next lemma, we estimate the degree of a consistent polynomial.
\begin{Lem}
\label{Lem2wn<1}
If the number $w\in (1/2;1)$ is consistent with a polynomial of degree $N\ges2$, then the next two inequalities are met:
$N >\log_w(1/2)$ (that is equivalent to $2w^N < 1$), and $N \ges \log_w(2-1/w)$.
\end{Lem}
{\bf Proof.}
1) Let us put $n = [- \log_2(-\log_2 w)]+1$.
Then $w\in[1/2^{1/2^{n-1}}, 1/2^{1/2^n})$.
So, according to item~3) of the theorem~\ref{TeorSqrtC}, 
the equality $N = 2^{n-1}(\deg F_d + 1) - 1$ is true, where $d=w^{2^{n-1}} < 1$.
Since $\deg F_d = 1$ only if $d=1$, then $\deg F_d \ges 2$.
From here 
$N \ges 2^{n-1}\cdot3 - 1 \ges 2^n > \log_w(1/2)$.
The first inequality is proved.

2) Let the polynomial $P(x) = c_0+c_1x+\ld+c_nx^N$ be consistent with the number $w$.
Then $P(w) = 0$ and the following relations are true:
$$
0 = P(w) \ges 1-w-w^2 - \ld-w^N = (1-2w+w^{N+1})/(1-w).
$$
From here 
$w^{N+1} \les 2w-1$,
so $N \ges \log_w(2w-1)-1 = \log_w(2-1/w)$.

That's what we needed to prove. 
\qed 


\subsection{\boldmath Negativity of the derivative of consistent functions and uniqueness of consistent points at the segment $[1/2;1]$}

First, we present a lemma containing the technical estimates, which is necessary to prove the theorem~\ref{TeorProizEdinSogl}.
\begin{Lem}
\label{LemTexOcenki}
1) For all $x\in[1/2; 0{.}6]$ the next inequality is true:
$-1+10x^4-8x^5 < 0$.

2) For all $x\in[1/2; 0{.}62]$ the next two inequalities are true:

$-1+8x^3-16x^4+20x^5-10x^6 < 0$,

$-1+8x^3-6x^4-12x^5+24x^6-12x^7 < 0$.

3) For all $x\in[0{.}61; 0{.}665]$ the next inequality is true:
$-1+6x^2-12x^3+6x^4+12x^5-10x^6 < 0$.

4) For all $x\in[0{.}61; 0{.}7]$ the next two inequalities are true:

$-1+6x^2-12x^3+16x^4-20x^5+10x^6+16x^7-14x^8 < 0$,

$-1+6x^2-12x^3+16x^4-20x^5+24x^6-28x^7+14x^8+20x^9-18x^{10} < 0$.

5) For all $x\in[0{.}61; 0{.}721]$ the next inequality is true:

$-1+6x^2-12x^3+16x^4-20x^5+24x^6-28x^7+32x^8-36x^9+40x^{10} - 20x^{11} < 0$.
\end{Lem}
{\bf Proof.} The derivative will be applied in the proof.

1) Let $x\in[1/2; 0{.}6]$. 
Put $y_1(x) = -1+10x^4-8x^5$.
Therefore, $y_1'(x) = 40x^3 (1-x) > 0$, so $y_1(x)$ increases.
So $y_1(x) \les y_1(0{.}6) = -0{.}32608 < 0$.

2) Let $x\in[1/2; 0{.}62]$.
Put $y_2(x) = -1+8x^3-16x^4+20x^5-10x^6$ and
$z_2(x) = -1+8x^3-6x^4-12x^5+24x^6-12x^7$.
Then $y_2'(x) = 4x^2(1-x)(15x^2-10x+6) > 0$ 
and
$$
z_2'(x) = x^2(1-x)(19{.}5+3(x-1/2)+66 (x-1/2)^2+84(x-1/2)^3) > 0.
$$
Therefore, the functions $y_2(x)$ and $z_2 (x)$ increase.
For this reason
$y_2(x) \les y_2(0{.}62) = -0{.}1933264518 < 0$ and
$z_2(x) \les z_2(0{.}62) = -0{.}1387036567 < 0$.

3) Let $x\in[0{.}61; 0{.}665]$ and $y_3(x) = -1+6x^2-12x^3+6x^4+12x^5-10x^6$. Then 
$y_3'(x) = x(1-x)(7{.}5+21(x-1/2)+90(x-1/2)^2+60(x-1/2)^3) > 0$,
so $y_3(x)$ is increasing.
Consequently $y_3(x) \les y_3(0{.}665) = -0{.}0064637084 < 0$.

4) Let $y_4(x) = -1+6x^2-12x^3+16x^4-20x^5+10x^6+16x^7-14x^8$ and
$z_4(x) = -1+6x^2-12x^3+16x^4-20x^5+24x^6-28x^7+14x^8+20x^9-18x^{10}$.
Then at $x\in[0{.}61; 0{.}7]$ the next relations are correct:
$y_4'(x) = x(1-x)(6+6(x-1/2)+90(x-1/2)^2+220(x-1/2)^3+280(x-1/2)^4+112(x-1/2)^5) > 0$ 
and
$z_4'(x) = x(1-x)( (6{.}8125-2{.}3125x) +54{.}125(x-1/2)^2+221{.}75(x-1/2)^3+
591{.}5(x-1/2)^4+833(x-1/2)^5+630(x-1/2)^6+180(x-1/2)^7 ) > 0$.
Therefore, the functions $y_4(x)$ and $z_4 (x)$ increase.
It follows that 
$y_4(x) \les y_4(0{.}7) = -0{.}00871334 < 0$ 
and
$z_4(x) \les z_4(0{.}7) = -0{.}0724555682 < 0$.

5) Let $x\in[0{.}61; 0{.}721]$. Put 
$y_5(x) = -1+6x^2-12x^3+16x^4-20x^5+24x^6-28x^7+32x^8-36x^9+40x^{10}-20x^{11}$.
Therefore 
$y_5'(x) = x(1-x)( (6{.}421875-0{.}9375x) +49{.}125(x-1/2)^2+179{.}25(x-1/2)^3+
519(x-1/2)^4+915(x-1/2)^5+1054(x-1/2)^6+700(x-1/2)^7+220(x-1/2)^8) > 0$.
So $y_5(x)$ is increasing.
Then $y_5(x) \les y_5(0{.}721) = -0{.}0517022846 < 0$.
\qed

\begin{Teor}[on derivative and uniqueness]
\label{TeorProizEdinSogl}
Let the function $F$ be consistent with the point $w$ lying in the segment $[1/2;1]$.
Then $F'(w) < 0$ and $w$ is the unique point in $[1/2;1]$ that is consistent with $F$.
\end{Teor} 
{\bf Proof.} 
If $w=1$, then obviously $F(x)=1-x$. In this case, $F'(w)=-1 < 0$.
Also, $w=1$ is the unique point consistent with the $1-x$ polynomial, since it has no other roots.

Let now $w\in[1/2;1)$. Denote $w_1=w$.
In order to prove the theorem, it is sufficient to show that if the point $w_2\in[w_1;1)$ is also consistent with $F$, then $w_1=w_2$ and $F'(w_1) < 0$ (both in the case when $F$ is a polynomial and in the case when $F$ is a series).

Since $1-w_1>0$, then $F(x)$ begins with the members $1-x-x^2\pm\ld$.
For further proof, we divide the set of all consistent functions $F$ with the specified beginning into three types.
Then for functions of the first and second types, we prove that $F'(x)<0$ for any $x$ in the segment $[w_1,w_2]$.
It follows that $F$ decreases on $[w_1,w_2]$.
Therefore, since $F(w_1)=F(w_2)=0$ (by virtue of the theorem~\ref{TeorFw(w)=0}), we get the necessary relations $w_1=w_2$ and $F'(w_1) < 0$.
For functions of the third type, we prove these relations using the results obtained for the first two types.

1) Let us assign {\bf the first type} to the functions $F$ of the form $F(x) = 1-x-x^2-x^3\pm\ld$.
Then $1-w_i-w_i^2>0$ for $i=1,2$.
Hence $w_1, w_2\in[1/2, (\sqrt5-1)/2)$ and the inclusions $[w_1,w_2]\sbs[1/2, (\sqrt5-1)/2)\sbs[1/2; 0{.}62]$ are true.

1.1) Consider the first subtype, where $F(x) = 1-x-x^2-x^3-x^4\pm\ld$. 
Then $1-w_i-w_i^2-w_i^3>0$ for $i=1,2$. 
Hence, $w_1,w_2\in[1/2; 0{.}6]$.
Next, evaluate the derivative $F'(x)$ for $x\in[w_1,w_2]\sbs[1/2; 0{.}6]$, using item~1) of lemma~\ref{LemTexOcenki}:
$$
F'(x) = -1-2x-3x^2-4x^3\pm\ld \les -1-2x-3x^2-4x^3 + \sum_{n=5}^\infty nx^{n-1} =
$$
$$ 
= (1+x+x^2+\ld)' - 2(1+2x+3x^2+4x^3) =
1/(1-x)^2-2(1+2x+3x^2+4x^3) =
$$
$$
= (-1+10x^4-8x^5)/(1-x)^2 < 0.
$$

1.2) Consider the second subtype, where $F(x) = 1-x-x^2-x^3+x^4\pm\ld$.

1.2.1) Consider a subtype, where $F(x) = 1-x-x^2-x^3+x^4-x^5\pm\ld$.
Again, evaluate the derivative $F'(x)$ for $x\in[w_1,w_2]\sbs[1/2; 0{.}62]$, using now item~2) of lemma~\ref{LemTexOcenki}:
$$
F'(x) = -1-2x-3x^2+4x^3-5x^4\pm\ld \les -1-2x-3x^2+4x^3-5x^4 + \sum_{n=6}^\infty nx^{n-1} =
$$
$$
= (1+x+x^2+\ld)' - 2(1+2x+3x^2+5x^4) =
1/(1-x)^2-2(1+2x+3x^2+5x^4) = 
$$
$$
= (-1+8x^3-16x^4+20x^5-10x^6)/(1-x)^2 < 0.
$$

1.2.2) Consider the last remaining subtype of the first type.
Let $F$ have the form $F(x) = 1-x-x^2-x^3+x^4+x^5\pm\ld$.
Due to the relations $1-w_i-w_i^2-w_i^3+w_i^4+w_i^5 = (1-w_i^3) (1-w_i-w_i^2) > 0$ when $i=1,2$, in fact $F$ has the form $F (x) = 1-x-x^2-x^3+x^4+x^5-x^6\pm\ld$. 
Evaluate $F'(x)$ for $x\in[w_1,w_2]\sbs[1/2; 0{.}62]$, using item~2) of lemma~\ref{LemTexOcenki}:
$$
F'(x) = -1-2x-3x^2+4x^3+5x^4-6x^5\pm\ld \les
$$
$$
\les-1-2x-3x^2+4x^3+5x^4-6x^5 + \sum_{n=7}^\infty nx^{n-1} =
$$
$$
= (1+x+x^2+\ld)' - 2(1+2x+3x^2+6x^5) = 1/(1-x)^2-2(1+2x+3x^2+6x^5) =
$$
$$
= (-1+8x^3-6x^4-12x^5+24x^6-12x^7)/(1-x)^2 < 0.
$$

2) Let us assign {\bf the second type} to the functions $F$ of the form $F(x) = 1-x-x^2+x^3\pm\ld$,
with some exceptions, which will be discussed below.
Then $1-w_i-w_i^2<0$ for $i=1,2$.
Hence, $w_1,w_2\in[(\sqrt5-1)/2; 1)$, so $[w_1, w_2]\sbs[(\sqrt5-1)/2; 1)\sbs[0{.}61; 1)$.
Since $1-w_i-w_i^2+w_i^3 = (1-w_i) (1-w_i^2) > 0$ for $i=1,2$, 
then in fact all functions of the second type have the form $F(x) = 1-x-x^2+x^3-x^4\pm\ld$. 

2.1) Consider a subtype where $F(x) = 1-x-x^2+x^3-x^4-x^5\pm\ld$.
Then $1-w_i-w_i^2+w_i^3-w_i^4>0$ for $i=1,2$, so $w_1, w_2\in[0{.}61; 0{.}665]$.
Evaluate $F'(x)$ for $x\in[w_1,w_2]\sbs[0{.}61; 0{.}665]$, using item~3) of lemma~\ref{LemTexOcenki}:
$$
F'(x) = -1-2x+3x^2-4x^3-5x^4\pm\ld \les
$$
$$
\les-1-2x+3x^2-4x^3-5x^4 + \sum_{n=6}^\infty nx^{n-1} =
$$
$$
= (1+x+x^2+\ld)' - 2(1+2x+4x^3+5x^4) = 1/(1-x)^2-2(1+2x+4x^3+5x^4) =
$$
$$
= (-1+6x^2-12x^3+6x^4+12x^5-10x^6)/(1-x)^2 < 0.
$$

2.2) Consider a subtype where $F(x) = 1-x-x^2+x^3-x^4+x^5\pm\ld$.
A~part of this subtype where $F$ has the form
\begin{equation} \label{eqIskl1}
F(x) = 1-x-x^2+x^3-x^4+x^5+x^6\pm\ld,
\end{equation}
we refer to the third type and consider below. 
So here we will only consider the case where $F$ has the form $F(x) = 1-x-x^2+x^3-x^4+x^5-x^6\pm\ld$.

2.2.1) Consider the subtype where $F(x) = 1-x-x^2+x^3-x^4+x^5-x^6-x^7\pm\ld$.
Then $1-w_i-w_i^2+w_i^3-w_i^4+w_i^5-w_i^6 = (1-w_i^2+w_i^3)(1-w_i-w_i^3) > 0$ for $i=1,2$. 
Hence, $w_1,w_2\in[0{.}61; 0{.}7]$.
Estimate $F'(x)$ for points $x\in[w_1,w_2]\sbs[0{.}61; 0{.}7]$, using item~4) of lemma~\ref{LemTexOcenki}:
$$
F'(x) = -1-2x+3x^2-4x^3+5x^4-6x^5-7x^6\pm\ld \les
$$
$$
\les-1-2x+3x^2-4x^3+5x^4-6x^5-7x^6 + \sum_{n=8}^\infty nx^{n-1} =
$$
$$
= (1+x+x^2+\ld)' - 2(1+2x+4x^3+6x^5+7x^6) =
$$
$$
= 1/(1-x)^2-2(1+2x+4x^3+6x^5+7x^6) =
$$
$$
= (-1+6x^2-12x^3+16x^4-20x^5+10x^6+16x^7-14x^8)/(1-x)^2 < 0.
$$

2.2.2) Consider the subtype where $F(x) = 1-x-x^2+x^3-x^4+x^5-x^6+x^7\pm\ld$.
A part of this subtype where $F$ has the form
\begin{equation} \label{eqIskl2}
F(x) = 1-x-x^2+x^3-x^4+x^5-x^6+x^7+x^8\pm\ld, 
\end{equation}
we also refer to the third type and consider below. 
So here we consider only the case where $F(x) = 1-x-x^2+x^3-x^4+x^5-x^6+x^7-x^8\pm\ld$.

2.2.2.1) Consider a subtype where the consistent function $F$ has the form 
$F(x) = 1-x-x^2+x^3-x^4+x^5-x^6+x^7-x^8-x^9\pm\ld$.
Then $1-w_i-w_i^2+w_i^3-w_i^4+w_i^5-w_i^6+w_i^7-w_i^8 > 0$ for $i=1,2$. 
Hence, $w_1,w_2\in[0{.}61; 0{.}7]$.
Evaluate $F'(x)$ for $x\in[w_1,w_2] \sbs[0{.}61; 0{.}7]$, using item~4) of lemma~\ref{LemTexOcenki}:
$$
F'(x) = -1-2x+3x^2-4x^3+5x^4-6x^5+7x^6-8x^7-9x^8\pm\ld \les
$$
$$
\les-1-2x+3x^2-4x^3+5x^4-6x^5+7x^6-8x^7-9x^8 + \sum_{n=10}^\infty nx^{n-1} =
$$
$$
= (1+x+x^2+\ld)' - 2(1+2x+4x^3+6x^5+8x^7+9x^8) =
$$
$$
= 1/(1-x)^2-2(1+2x+4x^3+6x^5+8x^7+9x^8) =
$$
$$
= (-1+6x^2-12x^3+16x^4-20x^5+24x^6-28x^7+14x^8+20x^9-18x^{10})/(1-x)^2 < 0.
$$

2.2.2.2) Consider the last subtype of the second type, where $F$ has the form $F(x) = 1-x-x^2+x^3-x^4+x^5-x^6+x^7-x^8+x^9\pm\ld$.
A part of this subtype where $F$ has the form
\begin{equation} \label{eqIskl3}
F(x) = 1-x-x^2+x^3-x^4+x^5-x^6+x^7-x^8+x^9+x^{10}\pm\ld,
\end{equation}
we will also refer to the third type.
So here we will analyze only the variant where $F$ has the form $F(x) = 1-x-x^2+x^3-x^4+x^5-x^6+x^7-x^8+x^9-x^{10}\pm\ld$.
Thus $1-w_i-w_i^2+w_i^3-w_i^4+w_i^5-w_i^6+w_i^7-w_i^8+w_i^9 = (1-w_i)(1-w_i^2-w_i^4-w_i^6-w_i^8) > 0$ for $i=1,2$. 
Hence, $w_1,w_2\in[0{.}61; 0{.}721]$.
Estimate $F'(x)$ for points $x\in[w_1,w_2] \sbs[0{.}61; 0{.}721]$, using item~5) of lemma~\ref{LemTexOcenki}:
$$
F'(x) = -1-2x+3x^2-4x^3+5x^4-6x^5+7x^6-8x^7+9x^8-10x^9\pm\ld \les
$$
$$
\les-1-2x+3x^2-4x^3+5x^4-6x^5+7x^6-8x^7+9x^8-10x^9 + \sum_{n=11}^\infty nx^{n-1} =
$$
$$
= (1+x+x^2+\ld)' - 2(1+2x+4x^3+6x^5+8x^7+10x^9) = 
$$
$$
= 1/(1-x)^2-2(1+2x+4x^3+6x^5+8x^7+10x^9) =
$$
$$
= (-1+6x^2-12x^3+16x^4-20x^5+24x^6-28x^7+
$$
$$
+32x^8-36x^9+40x^{10}-20x^{11})/(1-x)^2 < 0.
$$

So, if the function $F$ has the form of the first or (except for three cases) the second type, and the points $w_1\in[1/2;1)$, $w_2\in[w_1;1)$ are consistent with $F$, then $w_1=w_2$ and $F'(w_1)<0$.

3) Let us assign {\bf the third type} to the functions $F$, which were deferred when considering the second type, 
that is the functions having the form \eqref{eqIskl1}, \eqref{eqIskl2} or \eqref{eqIskl3}.

If $F$ has the form~\eqref{eqIskl1}, that is $F(x) = 1-x-x^2+x^3-x^4+x^5+x^6\ld$, 
then $1-w_i-w_i^2+w_i^3-w_i^4+w_i^5 = (1-w_i) (1-w_i^2-w_i^4) < 0$ for $i=1,2$. 
Hence, $w_1,w_2\in[0{.}71; 1)$.

If $F$ has the form~\eqref{eqIskl2}, that is, $1-x-x^2+x^3-x^4+x^5-x^6+x^7+x^8\ld$, 
then $1-w_i-w_i^2+w_i^3-w_i^4+w_i^5-w_i^6+w_i^7 = (1-w_i)(1-w_i^2-w_i^4-w_i^6) < 0$ for $i=1,2$. 
Hence, $w_1,w_2\in[0{.}71; 1)$.

If $F$ has the form~\eqref{eqIskl3}, that is, $F(x) = 1-x-x^2+x^3-x^4+x^5-x^6+x^7-x^8+x^9+x^{10}\ld$, 
then $1-w_i-w_i^2+w_i^3-w_i^4+w_i^5-w_i^6+w_i^7-w_i^8+w_i^9 = (1-w_i)(1-w_i^2-w_i^4-w_i^6-w_i^8) < 0$ for $i=1,2$. 
Hence, $w_1,w_2\in[0{.}71; 1)$.

So, in all the remaining variants we have $w_1,w_2\in[0{.}71; 1)\sbs[1/\sqrt2; 1)$.
Next, put $n = [- \log_2(-\log_2 w_1)]$.
Then $n\in\NN$, $d_1=w_1^{2^n}\in[1/2, 1/\sqrt2)$ and $d_2=w_2^{2^n}\in[d_1; 1)$.
Therefore, the function $F_{d_1}$, which is consistent with the point $d_1$, is of the first or second type.
Since $w_1, w_2\in[1/\sqrt[2^n]{2}, 1)$, then according to item~3) of the theorem~\ref{TeorSqrtC}, 
the equality is satisfied
\begin{equation} \label{eqFProdFd}
F(x) = \prod_{k=0}^{n-1} (1-x^{2^k}) \cdot F_{d_1}(x^{2^n}) = 
\prod_{k=0}^{n-1} (1-x^{2^k}) \cdot F_{d_2}(x^{2^n}).
\end{equation}
Therefore, the functions $F_{d_1}$ and $F_{d_2}$ are the same.
So, according to the proven for the first and second types, $d_1=d_2$ and $F_{d_1}'(d_1) < 0$.
So $w_1=w_2$.
It remains to show that $F'(w_1)<0$.
From the formula~\eqref{eqFProdFd} we get:
$$
F'(w_1) = \Bigl(\prod_{k=0}^{n-1} (1-x^{2^k})\Bigr)'\Bigr|_{x=w_1} F_{d_1}(d_1) +
\prod_{k=0}^{n-1} (1-w_1^{2^k}) \cdot F_{d_1}'(d_1) \cdot 2^n w_1^{2^n-1}.
$$
Hence, taking into account the relations $F_{d_1}(d_1)=0$ and $F_{d_1}'(d_1)<0$, we find that $F'(w_1) = \prod_{k=0}^{n-1} (1-w_1^{2^k}) \cdot F_{d_1}'(d_1) \cdot 2^n w_1^{2^n-1} < 0$.
The theorem is proved. 
\qed 


\subsection{Consistent polynomials for one particular sequence of numbers}

\begin{Teor} (cf.~with the theorem~\ref{TeorTabor} of J.~Tabor and J.~Tabor \cite[Theorem~3.1]{Tabor2}).
\label{TeorEdinPolozKor}

1) For each $n\in\NN$, the polynomial $P_n (x) = 1-x-x^2 - \ld-x^n$ has a single positive root $u_n$, 
moreover $u_1=1$ and $u_n\in (1/2;1/\sqrt{2})$ for any $n\ges 2$.

2) The sequence $\{u_n\}$ decreases and has asymptotics 
\begin{equation} \label{eqwn>0Asimp}
u_n = 1/2+1/2^{n+2}+\underline{O}(n/2^{2n}) \quad ({n\to\infty}).
\end{equation}

3) For each $n\in\NN$, the polynomial $P_n(x)$ is consistent with the point $u_n$.
\end{Teor}
{\bf Proof} is similar to the proof of item~4) of the theorem~\ref{TeorSoglW<1/2>1}:

1) The case $n=1$ is obvious. For any $x\ges 0$ we have: $P_n'(x) =-1-2x-\ld-nx^{n-1}<0$, so $P_n(x)$ decreases.
Also, if $n\ges 2$ then $P_n(1/\sqrt2) = (\sqrt2+1-\sqrt2^{n+1})/\sqrt2^{n} < 0$ and $P_n (1/2) = 1/2^n > 0$.
Therefore, for any $n\ges 2$ the polynomial $P_n$ has the only positive root $u_n$, and $u_n\in(1/2;1/\sqrt{2})$.

2) If $n\in\NN$ then $P_{n+1}(u_n) = P_{n}(u_n) - u_n^{n+1} = - u_n^{n+1} < 0 = P_{n+1}(u_{n+1})$.
Hence, due to the polynomial $P_{n+1}(x)$ decreases by $x$, the inequality $u_n > u_{n+1}$ follows for all $n\in\NN$.
So, the sequence $\{u_n\}$ decreases.

In order to obtain the asymptotic formula~\eqref{eqwn>0Asimp}, 
first let us find for any $a\in\RR$ the asymptotics of the sequence $\{Q_{n}(x_n)\}$, 
where $Q_n (x) = 1-2x +x^n$, $x_n = 1/2+1/2^{n+1}+an / 2^{2n}$ (up to the members of the order $n/2^{2n}$):
$$
Q_{n}(x_n) = 1 - 2(1/2+1/2^{n+1}+an/2^{2n}) + (1/2+1/2^{n+1}+an/2^{2n})^{n} =
$$
$$
= -1/2^n - an/2^{2n-1} - e^{n\ln(1+1/2^{n}+an/2^{2n-1})}/2^n =
(1-2an)/2^{2n} + \overline{o}(n/2^{2n}).
$$
So $P_{n}(x_{n+1}) = Q_{n+1}(x_{n+1})/(1-x_{n+1}) = (1-2a(n+1))/2^{2n+1} + \overline{o}(n/2^{2n})$.
From here, for sufficiently large values of $n$, we get the relations:

if $a=-1$, then 
$P_{n}(x_{n+1}) = (3+2n)/2^{2n+1} + \overline{o}(n/2^{2n}) > 0 = P_{n}(u_n)$,
so $u_n > x_{n+1} = 1/2+1/2^{n+2}-(n+1)/2^{2n+2}$;

if $a=1$, then 
$P_{n}(x_{n+1}) = (- 1-2n)/2^{2n+1} + \overline{o}(n/2^{2n}) < 0 = P_{n}(u_n)$,
so $u_n < x_{n+1} = 1/2+1/2^{n+2}+(n+1)/2^{2n+2}$.

Hence, the asymptotic equality~\eqref{eqwn>0Asimp} is satisfied.

3). Let $1\les k\les n$. Then due to the decreasing of the polynomial $P_{k-1}$ and the sequence $\{u_n\}$, we have:
$c_k\cdot(c_0 + c_1 u_n +\ld+ c_{k-1} u_n^{k-1}) = - P_{k-1}(u_n) < - P_{k-1}(u_{k-1}) = 0$.
So, by definition~\ref{OprSoglFunc}, the polynomial $P_n(x)$ is consistent with the point $u_n$.
That's what we needed to prove. 
\qed

\subsection{\boldmath The set of points which are consistent with polynomials, is everywhere dens on the segment $[1/2;1]$}

\begin{Teor}
\label{TeorPlotnSoglMnog}
On the segment $[1/2; 1]$, the set of points which are consistent with polynomials, is countable and everywhere dense.
\end{Teor}
{\bf Proof.}
Let us denote the set of points $w\in[1/2; 1]$ that are consistent with polynomials, by $W_P$.

1) First, we show that the set $W_P$ is countable.

By definition~\ref{OprSoglFunc}, any point that is consistent with some polynomial, is its root.
Since the coefficients of consistent polynomials are $\pm1$, the set of such polynomials is no more than countable. 
Therefore, due to each consistent polynomial has only finite number of roots, the set $W_P$ is no more than countable too.

On the other hand, by the theorem~\ref{TeorEdinPolozKor}, the polynomial $1-x-x^2 - \ld-x^n$ for any $n\in\NN$ has the root $u_n\in(1/2;1/\sqrt{2})$, and $u_n$ is consistent with this polynomial.
Therefore, the set $W_P$ is at least countable.

2) In order to prove the density of the set $W_P$ on the segment $[1/2; 1]$, 
we just need to show that for any $\eps>0$ any subinterval $(w-\eps, w+\eps)\sbs(1/2, 1)$ contains some point,
which is consistent with a polynomial. 
In the case, where the point $w$ itself is consistent with the polynomial, $w$ is the desired point.
So consider the remaining case where $w$ is consistent with some series $F(x) = c_0+c_1x+c_2x^2+\ld$.

According to proposition~\ref{PredlAnalF}, 
the convergence region of the series $F(z) = c_0+c_1z+c_2z^2+\ld$ is the circle $\{z\in\CC \bigm| |z|<1\}$.
So due to the properties of power series (see~\cite[Ch.~II, \S3, p.~6, p.~66]{Privalov}), 
the polynomials $F_n(z) = c_0+c_1z+\ld+c_nz^n$, $n=0,1,2,\ld$ converge to the function $F(z)$ 
uniformly in the circle $\{z\in\CC \bigm| |z| \les 1-\eps/2\}$, and thus in circle $K_\eps=\{z\in\CC \bigm| |z-w|\les\eps/2\}$.

Next, $F(w)=0$ by the theorem~\ref{TeorFw(w)=0}.
Since $F(z)$ is not identically zero, then, according to the uniqueness theorem 
(see~\cite[Ch.~V, \S2, p.~4, p.~202]{Privalov}), the point $w$ is an isolated root of the function $F(z)$.
Therefore, we can consider the number $\eps>0$ so small that $F(z)$ in the circle $K_\eps$ has no roots other than $w$.
This means that the value $\mu_\eps =\inf_{|z-w|=\eps/2} |F(z)|$ is positive.
Due to the uniform convergence of $F_n (z)$ to $F(z)$ in the circle $K_\eps$, 
there exists such a number $N\in\NN$ that the inequality $|F_N(z)-F(z)| < \mu_\eps$ 
is satisfied on the contour $C_\eps=\{z\in\CC \mid |z-w|=\eps/2\}$.

By virtue of the theorem~\ref{TeorProizEdinSogl}, we have $F'(w)<0$, so $w$ is the root of $F(z)$ of multiplicity~$1$.
Now let us apply Rouchet's theorem (see~\cite[Ch.~VII, \S2, p.~2, p.~245]{Privalov}).
Since $|F_N(z)-F(z)| <\mu_\eps \les |F(z)|$ on the contour $C_\eps$, 
then the functions $F(z)$ and $F (z)+(F_N(z)-F (z))=F_N (z)$ have the same number of zeros (taking into account the multiplicity)
inside $C_\eps$ (that is, inside the circle $K_\eps$).
Since $z=w$ is the only root (of multiplicity $1$) of $F(z)$ inside the circle $K_\eps$, 
then $F_N(z)$ also has exactly one root $z_N$ inside $K_\eps$, and its multiplicity equals $1$.

The unitary polynomial $F_N(z)$ has real coefficients, 
so the number $\overline{z_N}$, which is complex conjugate to $z_N$, is also the root of $F_N(z)$. 
It is clear that $\overline{z_N}$, like $z_N$, lies inside the circle $K_\eps$, where $F_N(z)$ has only one root. 
So $\overline{z_N}=z_N$, that is, $z_N$ is real.
Hence, $z_N \in (w-\eps, w+\eps)$.

Now let us find in the interval $(w-\eps, w+\eps)$ a point that is consistent with some polynomial 
(we do not claim that this is the point $z_N$).
To do this, we form the set $Z_\eps$ of all those roots of the polynomials $F_0,\ld,F_N$, 
which fall within the interval $(w-\eps, w+\eps)$. 
$Z_\eps$ contains, in particular, the point $z_N$.
In the set $Z_\eps$, take the element $z_*$ closest to the point $w$, 
and then select from $F_0,\ld,F_N$ the polynomial $F_m$ of minimal degree $m$, which has $z_*$ as its root.
Let us show that the point $z_*$ is consistent with the polynomial $F_m$ (and is the desired one).
Indeed, according to the construction of $z_*$ and $F_m$, 
the polynomials $F_0,\ld,F_{m-1}$ have no roots on the segment between the points $z_*$ and $w$. 
Therefore, the values $F_0(z_*),\ld, F_{m-1}(z_*)$ have the same signs as $F_0(w),\ld,F_{m-1}(w)$. 
Further, since the point $w$ is consistent with the series $F$, then by definition~\ref{OprSoglFunc} for any $k=1,\ld, m$ inequalities~\eqref{eqSoglKoef} are met, that is
$$
c_k\cdot(c_0 + c_1 w +\ld+ c_{k-1} w^{k-1}) = c_k F_{k-1}(w) < 0.
$$
Since for all $k=1,\ld,m$ the values $F_{k-1}(w)$ and $F_{k-1}(z_*)$ have the same signs, the inequalities follow 
$$
c_k F_{k-1}(z_*) = c_k\cdot(c_0 + c_1 z_* +\ld+ c_{k-1} z_*^{k-1}) < 0.
$$
Therefore, by definition~\ref{OprSoglFunc}, the polynomial $F_m$ is consistent with the point $z_*\in(w-\eps, w+\eps)\sbs(1/2, 1)$.
That's what we needed to prove. \qed

\subsection{\boldmath Single-ended coefficient-wise limits of functions that are consistent with points from the segment $[1/2;1]$}

In this subsection, we describe the behavior of functions $F_w$ that are consistent with points $w\in[1/2;1]$, 
and also the behavior of attached to $F_w$ rows $F_w^{\pm}$ (see~definition~\ref{OprSoglFunc}), 
when $w$ tends to the limit on one side.

\begin{Teor}
\label{TeorPredelSoglFunc}
1) For any $w_0\in[1/2; 1)$ in the sense of coefficient-wise convergence, 
the next chain of equalities  is true: 
$\lim_{w\to w_0 \pm 0} F_w = \lim_{w\to w_0 \pm 0} F_w^- = \lim_{w\to w_0 \pm 0} F_w^+ = F_{w_0}^{\pm}$, 
where the function $F_w$ is consistent with $w$, and the series $F_w^+$, $F_w^-$ are attached to $F_w$ .
Moreover, if the point $w_0$ is consistent with a series (not with a polynomial), 
then $\lim_{w\to w_0} F_w = \lim_{w\to w_0} F_w^{\pm} = F_{w_0}$.

2) The following formula is сorrect:
$\lim_{w\to 1-0} F_w(x) = \lim_{w\to 1-0} F_w^{\pm}(x) = \prod_{k=0}^{\infty}(1-x^{2^k})$.
\end{Teor}
{\bf Proof.} 
1) Let ${w_0}\in[1/2;1)$.
For each $n\in\NN$, we can select the number $\de_n>0$ so that 
any point $w\in({w_0}-\de_n, {w_0}+\de_n)\setminus\{{w_0}\}$ is consistent either with a polynomial of degree at least $n$, or with a series.
This is possible because there are only a finite number of points that are consistent with polynomials of degree lower than $n$.
Then for any $w\in({w_0}-\de_n, {w_0}+\de_n)\setminus\{{w_0}\}$, the sets of the first $n$ coefficients of the functions $F_w$, $F_w^-$ and $F_w^+$ will be the same, as follows from the definition~\ref{OprSoglFunc}.
Therefore, the following limits are equal (if they exist): 
$\lim_{w\to {w_0} \pm 0} F_w = \lim_{w\to {w_0} \pm 0} F_w^- = \lim_{w\to {w_0} \pm 0} F_w^+$.
So all we have to do is to prove that $\lim_{w\to {w_0} \pm 0} F_w = F_{w_0}^{\pm}$.

1a) First, we prove the equality $\lim_{w\to {w_0}} F_w = F_{w_0}$ in the case when the number ${w_0}$ is consistent 
with some series $F_{w_0}(x) = c_0+c_1x+c_2x^2+\ld$.
For this purpose, it is sufficient for each $n\in\NN$ to find such a $\de\in(0;\de_n)$ that 
for any $w\in({w_0}-\de, {w_0}+\de)$ the first $n$ coefficients 
of the consistent functions $F_w$ and $F_{w_0}$ are the same.
By the definition~\ref{OprSoglFunc}, the required matching of coefficients is equivalent to performing inequalities
\begin{equation} \label{eqSoglKoef_1}
c_{k}(c_0+c_1w+\ld+c_{k-1}w^{k-1}) < 0, \quad k=0,1,\ld, n-1
\end{equation}
for all $w\in({w_0}-\de, {w_0}+\de)$.
Since the series $F_{w_0}$ is consistent with the point ${w_0}$, the inequalities~\eqref{eqSoglKoef_1} are satisfied when $w={w_0}$.
Due to the continuity of polynomials $c_{k}(c_0+c_1w+\ld+c_{k-1}w^{k-1})$ by $w$ , these inequalities will also be satisfied for all $w\in({w_0}-\de, {w_0}+\de)$ if $\de\in(0;\de_n)$ is small enough.

1b) Now we prove that if the number ${w_0}$ is consistent with the polynomial $P_{w_0}(x) = c_0+c_1x+\ld+c_Nx^N$, then $\lim_{w\to {w_0}-0} F_w = F_{w_0}^-$
(the proof of the equality $\lim_{w\to {w_0}+0} F_w = F_{w_0}^+$ is similar).
According to the definition~\ref{OprSoglFunc} we have:
$$
F_{w_0}^-(x) = c_0+c_1 x+\ld+c_Nx^N -c_0x^{N+1}-c_1x^{N+2}-\ld-c_Nx^{2N+1} - \ld = \\
$$
$$
= P_{w_0}(x)(1-x^{N+1}-x^{2(N+1)}-\ld) =\sum_{k=0}^\infty b_k x^k,       
$$
where $b_k = c_k$ when $k=0,\ld, N$;
$b_{i(N+1)+j} = - c_j$ when $i\in\NN$ and $j=0,\ld, N$.

It is enough for us to prove for any $n\in\NN$ the existence of such a $\de\in(0;\de_n)$ that 
for any $w\in({w_0}-\de, {w_0})$ the first $n$ coefficients of the consistent function $F_w$ and of the series $F_{w_0}^-$ coincide.
By the definition~\ref{OprSoglFunc}, the required matching of coefficients is equivalent to performing inequalities
\begin{equation} \label{eqSoglBi}
b_k(b_0+b_1w+\ld+b_{k-1}w^{k-1}) < 0
\end{equation}
for all $k=0,1,\ld, n-1$ and $w\in({w_0}-\de, {w_0})$.
Consider 3 cases:

I) If $0\leqslant k\leqslant N$, then the inequalities~\eqref{eqSoglBi} coincide with the inequalities~\eqref{eqSoglKoef_1} from item~1), and the proof is the same as there.

II) Let $k = i(N+1)$, where $i\in\NN$.
Then $b_k = -c_0 = -1$, and the inequality~\eqref{eqSoglBi} is equivalent to the inequality
$$
-(c_0+c_1w+\ld+c_Nw^N)(1-\sum_{p=1}^{i-1} w^{p(N+1)}) =
-P_{{w_0}}(w)(1 - \sum_{p=1}^{i-1} w^{p(N+1)}) < 0.
$$
We have: $P_{{w_0}}({w_0}) = 0$ by definition~\ref{OprSoglFunc}, $(P_{{w_0}})'({w_0}) < 0$ by the theorem~\ref{TeorProizEdinSogl} and $2a^N<1$ by the Lemma~\ref{Lem2wn<1}.
Therefore, there exist such a $\de\in(0;\de_n)$ that for all $w\in({w_0}-\de, {w_0})$ there will be $P_{{w_0}} (w) > 0$ and $2w^N<1$.
This leads to the necessary inequality:
$$
-P_{{w_0}}(w)(1-\sum_{p=1}^{i-1} w^{p(N+1)}) < - P_{{w_0}}(w)(1 - \sum_{p=1}^\infty w^{p(N+1)}) =
$$
$$
= - P_{{w_0}}(w)(1-2w^{N+1})(1-w^{N+1}) < 0.
$$

III) Let $k = i(N+1)+j$, where $i\in\NN$ and $1\leqslant j\leqslant N$.
Then
$$
b_k(b_0+b_1w+\ld+b_{k-1}w^{k-1}) =
$$
$$
= - c_j\Bigl( P_{{w_0}}(w)(1-\sum_{p=1}^{i-1} w^{p(N+1)}) - w^{i(N+1)}(c_0+c_1w+\ld+c_{j-1}w^{j-1}) \Bigr).
$$
Hence, for $w={w_0}$, due to the consistency of the point ${w_0}$ with $P_{{w_0}}$, we get:
$$
b_k(b_0+b_1w+\ld+b_{k-1}w^{k-1}) = {w_0}^{i(N+1)}c_j(c_0+c_1a+\ld+c_{j-1}{w_0}^{j-1}) < 0.
$$
In view of the continuity of $b_k(b_0+b_1w+\ld+b_{k-1}w^{k-1})$ by $w$, 
this inequality will also hold for all $w\in({w_0}-\de, {w_0})$ if $\de\in(0;\de_n)$ is small enough.

2) By virtue of item~3) of the theorem~\ref{TeorSqrtC}, $F_w (x) = \prod_{k=0}^{n-1}(1-x^{2^k}) F_d(x^{2^n})$ 
when $w\in[1/\sqrt[2^n]{2}, 1)$, where $d = w^{2^n}$.
Therefore, since $F_d(x^{2^n}) = 1 + c_1 x^{2^n} + \ldots$, it follows that for any $w\in[1/\sqrt[2^n]{2},1)$ 
the first $2^n$ coefficients of the series $F_w(x)$ and $\prod_{k=0}^{\infty}(1-x^{2^k}) $ coincide.
So $\lim_{w\to 1-0} F_w(x) = \prod_{k=0}^{\infty}(1-x^{2^k})$.
The equality of the limits $\lim_{w\to 1-0} F_w(x)^{\pm} = \lim_{w\to 1-0} F_w(x)$ can be proved in the same way as in point 1).
\qed

\subsection{\boldmath On lexicographic increase of consistent functions $F_w$ relative to $w\in[1/2;1)$}

In this subsection, we talk about lexicographically increasing dependence of the consistent functions $F_w$ for points $w\in[1/2;1)$.
\begin{Teor}
\label{TeorMonotSoglFunc}
If $1/2\leqslant a < b < 1$, then $F_a^+ \prec F_b^-$ (see~definition~\ref{OprLeksSravnShod}).
\end{Teor}
{\bf Proof.}
For any $i=0,1,2,\ld$ we will denote through $c_i(w)$ the coefficients of functions $F_w$, 
that is consistent with points $w\in\RR$, and through $c_i^\pm(w)$ the coefficients of attached functions $F_w^\pm$.

a) First, we prove the non-strict inequality $F_a^+ \preceq F_b^-$ from the opposite.
Let $F_a^+ \succ F_b^-$.
Then there exists such a number $n$ that $c_i^+(a)=c_i^-(b)$ when $i=0,1,\ld, n-1$, 
but $c_n^+(a) > c_n^-(b)$, that is $c_n^+(a)=1$ and $c_n^-(b)=-1$. 
According to the item~1) of the theorem~\ref{TeorPredelSoglFunc}, there is some $\de\in(0,(b-a)/2)$,
such that whatever $i=0,1,\ld,n$ is, for all $w\in(a,a+\de)$ the equality $c_i(w)=c_i^+(a)$ is true 
and for all $w\in(b-\de,b)$ another equality $c_i(w)=c_i^-(b)$ is true.
Let's select the points $w_1\in(a, a+\de)$ and $w_2\in(b-\de, b)$, 
which are consistent with the series (this is possible due to the theorem~\ref{TeorPlotnSoglMnog}).
Then $c_i(w_1)=c_i(w_2)$, $i=0,1,\ld, n-1$, but $c_n (w_1) > c_n(w_2)$.

Denote by $k$ the smallest number for which the coefficient $c_k (w)$ changes when $w$ increases from $w_1$ to $w_2$.
It is clear that $2\les k \les n$ (see~item~1 of the remark~\ref{ZamSoglFun}).
Then for $i=0,1,\ld, k-1$, the coefficients $c_i (w)=c_i$ do not depend on $w\in[w_1,w_2]$.
Put $P_{k-1}(x) = c_0+c_1x+\ld+c_{k-1}x^{k-1}$.
By virtue of the definition~\ref{OprSoglFunc}, we have the equalities 
$\sign(P_{k-1}(w_1)) = - c_k(w_1) = -1$ and $\sign(P_{k-1}(w_2)) = - c_k(w_2) = 1$.
Therefore, the polynomial $P_{k-1}$ has at least one root on $[w_1,w_2]$.
Then by definition~\ref{OprSoglFunc} any such root $d\in[1/2;1)$ is consistent with $P_{k-1}$.
By the theorem~\ref{TeorProizEdinSogl}, the consistent point $d$ is unique and $P_{k-1}'(d) < 0$.
Consequently, $P_{k-1}(w) > 0$ for any $w\in[w_1, d)$, and $P_{k-1}(w) < 0$ for any $w\in(d, w_2]$.
So $c_k(w_1) = - \sign(P_{k-1}(w_1)) = -1$ and $c_k(w_2) = -\sign(P_{k-1}(w_2)) = 1$.
Hence $k=n$, since $c_k(w_1)=c_k(w_2)$ when $k<n$.
Then the equalities $c_n(w_1)=-1$ and $c_n(w_2)=1$ are true, 
which contradicts the inequality $c_n(w_1) > c_n(w_2)$.
This contradiction implies that $F_a^+ \preceq F_b^-$.

b) Next, we prove that $F_a^+ \prec F_b^-$.
As in item~a), select the points $w_1<w_2$ in the interval $(a,b)$, 
that are consistent with the series $F_{w_1}$ and $F_{w_2}$ respectively.
Since each of these series is consistent with only one point (according to the theorem~\ref{TeorProizEdinSogl}), 
then $F_{w_1} \neq F_{w_2}$.
So the relations $F_a^+ \preceq F_{w_1} \prec F_{w_2} \preceq F_b^-$ are fulfilled 
according to the non-strict inequalities of item~a).
Hence, $F_a^+ \prec F_b^-$.
\qed

\begin{Zam}
At the point $w=1$, the lexicographic increment of consistent and attached functions, which is indicated by the theorem~\ref{TeorMonotSoglFunc}, is violated,
and these functions seem to make a jump down.
Indeed, by item~2) of the theorem~\ref{TeorPredelSoglFunc} we have:
$\lim_{w\to 1-0} F_w (x) = \prod_{k=0}^{\infty}(1-x^{2^k}) = 1-x-x^2+x^3-x^4+x^5+x^6-\ld$.
At the same time, since the point $1$ is consistent with the polynomial $1-x$, by virtue of the definition~\ref{OprSoglFunc} and the example~\ref{PrimSoglW1/sqrt2} we have:
$F_1^-(x) = 1-x-x^2+x^3-x^4+x^5-x^6+\ld = F_{1/\sqrt2}(x) \prec \lim_{w\to 1-0} F_w(x)$.
\end{Zam}


\subsection{\boldmath Several examples of consistent series and polynomials for points of the interval $(1/2;1)$}

\begin{Prim}
The point $(\sqrt5-1)/2\approx 0{.}618034$ is consistent with the unitary polynomial $P_2(x)=1-x-x^2$.
\end{Prim}
{\bf Proof.}
This is a special case of the theorem~\ref{TeorEdinPolozKor} for $n=2$.
\qed

\begin{Prim}
\label{PrimSoglW1/sqrt2}
The point $1/\sqrt2 \approx 0{.}707107$ is consistent with the unitary series 
$F_{1/\sqrt2}(x) = 1-x-x^2+x^3-x^4+x^5-x^6+\ld$, whose signs alternate starting with the term $-x^2$.
For $|x|<1$, its sum is equal to $F_{1/\sqrt2}(x) = (1-x)(1-2x^2)/(1-x^2) = (1-2x^2)/(1+x)$.
\end{Prim}
{\bf Proof.}
This follows from item~3) of the theorem~\ref{TeorSoglW<1/2>1} and item~2) of the theorem~\ref{TeorSqrtC}.
\qed

\begin{Prim}
\label{PrimSoglW_sqrt4_2}
The point $w = 1/\sqrt[4]{2} = \sqrt{1/\sqrt2} \approx 0.840896$ is consistent with the series 
$$
F_w(x) = (1-x)F_{1/\sqrt2}(x^2) = 1-x-x^2+x^3-x^4+x^5+x^6-x^7-x^8+x^9+x^{10}-\ldots.
$$
Here, the periodicity of signs begins with the term $-x^4$. 
For $|x|<1$, the sum of the series is equal to
$F_w(x) = (1-x)F_{1/\sqrt2}(x^2) = (1-x)(1-2x^4)/(1+x^2)$.
\end{Prim}
{\bf Proof.}
This follows from the example~\ref{PrimSoglW1/sqrt2} and item~2) of the theorem~\ref{TeorSqrtC}.
\qed

\begin{Prim}
\label{PrimSoglW2/3}
The number $2/3$ is consistent with some series $F_{2/3}(x)$ (not a polynomial). 
As the calculations show, its beginning has the form:
$$
F_{2/3}(x) = 1-x-x^2+x^3-x^4+x^5-x^6-x^7+x^8-x^9+x^{10}-x^{11}+x^{12}+
$$
$$
+x^{13}-x^{14}-x^{15}+x^{16}-x^{17}+x^{18}-x^{19}+x^{20}+x^{21}-x^{22}-x^{23}+x^{24} - \ldots.
$$
We were not able to find any regularity in the signs of the members of this series.
\end{Prim}
{\bf Proof.}
Consistency with a series, not a polynomial, follows from the proposition~\ref{PredlSoglIrraz}, 
since the number $2/3$ is rational and is not equal to $\pm1$.
\qed

\begin{Prim}
\label{PrimSoglW1/sqrt3}
The number $w = 1/\sqrt{3}\approx 0.577350$ is also consistent with a series. 
The beginning of this, according to the calculations, has the form:
$$
F_{1/\sqrt{3}}(x) = 1-x-x^2-x^3+x^4-x^5+x^6+x^7-x^8+x^9+x^{10}-x^{11}+x^{12}-
$$
$$
-x^{13}+x^{14}+x^{15}-x^{16}+x^{17}-x^{18}+x^{19}+x^{20}-x^{21}+x^{22}+x^{23}-x^{24} - \ldots.
$$
Patterns of signs were also not found here.
Here we did not find any regularity too.
\end{Prim}
{\bf Proof.}
If $w = 1/\sqrt{3}$ were consistent with the polynomial $c_0 + c_1 x +\ld+ c_n x^n$, 
then the rational number $w^2=1/3\neq\pm1$ would be the root of the unitary polynomial $c_0 + c_2x + c_4x^2 +\ld+ c_{2[n / 2]}x^{[n/2]}$.
But, according to the proposition~\ref{PredlSoglIrraz}, this is impossible.
\qed

\section{Relation between consistent and anti-consistent functions, and global extremes of functions in exponential Takagi class}
\label{SectSvazSoglSEkst}

The only theorem~\ref{TeorPoiskGlobExtr} of this section shows 
how consistent and anti-consistent polynomials and series can be used 
to find global extremes of functions $T_v$ in exponential Takagi class (see formula~\eqref{eqDefTv}). 

In order to prove the theorem~\ref{TeorPoiskGlobExtr}, we need the following lemma.
\begin{Lem}
\label{LemSvnLin}
Let $v\in(-1;1)$, $m\in\ZZ$, $n\in\{0,1,2,\ld\}$ and $[a, b]$ is one of two segments $[m/2^n-1/2^{n+1}, m / 2^n]$ or $[m/2^n, m / 2^n+1/2^{n+1}]$.
Then

a) the function $S_{v, n}(x) = \sum\limits_{k=0}^n v^k T_0(2^k x)$ is linear on $[a, b]$, and its derivative $S'_{v, n}$ is constant on the interval $(a,b)$;

b) the sets of points of the global maximum of the functions $S_{v, n}$ and $T_v$ on the segment $[a,b]$

- are in the right half of the segment $[a, b]$ if $S'_{v, n} > 0$;

- are in the left half of the segment $[a, b]$ if $S'_{v, n} < 0$;

- are symmetric about the middle of the segment $[a,b]$ if $S'_{v, n} = 0$.
\end{Lem}
{\bf Proof.}
It follows from the equality $T_0(x) = \rho(x,\ZZ)$, that 
the function $T_0$ is linear on the intervals $[m-1/2,m]$ and $[m,m+1/2]$ and is symmetric about the point $m/2$ 
(i.e. $T_0(m/2+x)=T_0(m/2-x)$ when ${x\in\RR}$).
It follows, that for any $n=0,1,2,\ld$ the function $T_0(2^n x)$ 
is linear on the segments $[m/2^n-1/2^{n+1}, m/2^n]$ and $[m/2^n, m/2^n+1/2^{n+1}]$,
and is symmetric about the points $m/2^{n+1}$.
Therefore, for any $m\in\ZZ$ and $n=0,1,2,\ld$ the function $S_{v,n}$ is linear 
on the segments $[m/2^n-1/2^{n+1}, m/2^n]$ and $[m/2^n, m/2^n+1/2^{n+1}]$, 
and the function $R_{v, n}(x) = \sum_{k=n+1}^\infty v^k T_0(2^k x)$ is symmetric about points $m/2^{n+2}$, 
including about the midpoints of the specified segments.
Hence, by virtue of the equality $T_v(x) = S_{v,n}(x)+R_{v,n}(x)$, the lemma statement follows.
\qed

\begin{Teor}
\label{TeorPoiskGlobExtr}
Let $v\in(-1;1)$ and $E_v$ be the set of points of the global maximum (respectively minimum) of the function $T_v$ on the segment $[0;1]$. Then the following statements are true:

1) If the point $2v$ is consistent (respectively anti-consistent) with a series 
$F_{2v}(x) = c_0 + c_1 x + \ldots + c_n x^n + \ldots$, then

1a) The set $E_v$ contains only two (possibly coinciding) points: 
$x^-(v)\in[0; 1/2]$ and $x^+(v)\in[1/2;1]$.
The first point and its binary expansion have the form
\begin{equation} \label{eqx^-(v)}
x^-(v) = 1/2 - F_{2v}(1/2)/4 = 0{.}x_1^-x_2^- \ld,
\end{equation}
where 
\begin{equation} \label{eqx_n^-}
x_n^- = (1-c_{n-1})/2, \quad n\in\NN.
\end{equation}
The second point and its binary expansion have the form
\begin{equation} \label{eqx^+(v)}
x^+(v) = 1 - x^-(v) = 1/2 + F_{2v}(1/2)/4 = 0{.}x_1^+x_2^+\ld, 
\end{equation}
where 
\begin{equation} \label{eqx_n^+}
x_n^+ = 1-x_n^- = (1+c_{n-1})/2, \quad n\in\NN.
\end{equation}

1b) The global maximum (respectively minimum) of the function $T_v$ can be calculated using the following formulae:
\begin{equation} \label{eqMaxT_v1}
T_v(x^\pm(v)) = \frac{1}{2(1-v)} - \frac{1}{4}\sum_{n=0}^\infty c_n\cdot(2v)^n\sum_{p=n}^\infty \frac{c_p}{2^p}
\end{equation}
and 
\begin{equation} \label{eqMaxT_v2}
T_v(x^\pm(v)) = \frac{1}{2(1-v)} - \frac{1}{4\pi i}\int_{|z|=r} \frac{F_{2v}(z)F_{2v}(v/z)}{2z-1}dz,
\end{equation}
where $r$ is any number from the interval $(\max (1/2, v), 1)$.

2) If the point $2v$ is consistent (respectively anti-consistent) with a polynomial
$P_{2v,N}(x) = c_0 + c_1 x + \ldots + c_n x^N$,
and $F_{2v}^\pm$ are attached series, then the following statements are true : 

2a) The set $E_{v,N}$ of points of the global maximum (respectively minimum) 
of the function $S_{v, N}(x) = \sum_{n=0}^N v^n T_0(2^nx)$ on the segment $[0;1]$ 
has the form $E_{v, N} = [a_N, b_N]\cup[1-b_N, 1-a_N]$,
where
\begin{equation} \label{eqabN}
a_N(v) = 1/2 - P_{2v, N}(1/2)/4 - 1/2^{N+2}, \quad
b_N(v) = a_N (v) + 1/2^{N+1}.
\end{equation}
In this case, the binary expansions of both points $a_N (v)$ and $1-a_N (v)$ are finite:
$a_N(v) = 0{.}x_1^-x_2^- \ld x_{N+1}^-$,
$1-a_N(v) = 0{.}x_1^+x_2^+ \ld x_{N+1}^+$,
where 
\begin{equation} \label{eqx_n^-Pol}
x_n^- = (1-c_{n-1})/2, \quad x_n^+ = 1-x_n^- = (1+c_{n-1})/2, \quad n=1,\ld, N+1.
\end{equation}
The global maximum (respectively minimum) $M_{v, N}$ of the function $S_{v, N}$ on the segment $[0;1]$ 
can be calculated using the formulae
\begin{equation} \label{eqMvN1}
M_{v, N} = \frac{1-v^{N+1}}{2(1-v)} - \frac{1}{4}\sum_{n=0}^N c_n\cdot (2v)^n \sum_{i=n}^N \frac{c_i}{2^i}
\end{equation}
and 
\begin{equation} \label{eqMvN2}
M_{v, N} = \frac{1-v^{N+1}}{2(1-v)} - \frac{1}{4\pi i}\int_{|z|=r} \frac{P_{v, N}(z)P_{v, N}(v/z)}{2z-1} dz,
\end{equation}
where $r$ is any number from the set $(0; 1/2)\cup (1/2;1)$.

2b) In the case of $v^{N+1}>0$ (that is, $v>0$ or $N$ is odd) the set $E_v$ is infinite.
Moreover, $E_v$ has Hausdorff dimension $1/(N+1)$ and consists of all the points $x$, 
that have the binary expansion of the form 
\begin{equation} \label{eqxExtrSoglPol1}
x = 0{.}
\left[\begin{array}{l}x_1^-x_2^-\ld x_{N+1}^-\\x_1^+x_2^+\ld x_{N+1}^+\end{array}\right.
\left[\begin{array}{l}x_1^-x_2^-\ld x_{N+1}^-\\x_1^+x_2^+\ld x_{N+1}^+\end{array}\right.
\ld.
\end{equation}
Any point $x$ of the set $E_v$ can also be written in the form
\begin{equation} \label{eqxExtrSoglPol3}
x = 1/2 \pm F(1/2)/4,
\end{equation}
where $F$ is some intermediate series for the number $2v$. 
And vice versa, any point of the form~\eqref{eqxExtrSoglPol3} belongs to $E_v$.
Besides that, the following formulae are correct:
\begin{equation} \label{eqinfEv}
\begin{array}{rl}
\inf E_v &= 0{.}x_1^-\ld x_{N+1}^-x_1^-\ld x_{N+1}^-\ld =\\
&= 1/2 - P_{2v, N}(1/2)/(4-1/2^{N-1}) = 1/2-F_{2v}^+(1/2)/4,
\end{array}
\end{equation}
\begin{equation} \label{eqinfEv012}
\begin{array}{l}
\sup(E_v\cap[0;1/2]) = 0{.}x_1^-\ld x_{N+1}^- x_1^+\ld x_{N+1}^+x_1^+\ld x_{N+1}^+\ld = \\
= 1/2-P_{2v, N}(1/2)\cdot(1-1/2^N)/(4-1/2^{N-1}) = 1/2 - F_{2v}^-(1/2)/4,
\end{array}
\end{equation}
\begin{equation} \label{eqsupEv}
\begin{array}{rl}
\sup E_v &= 1 - \inf E_v = 0{.}x_1^+\ld x_{N+1}^+x_1^+\ld x_{N+1}^+\ld =\\
&= 1/2 + P_{2v,N}(1/2)/(4-1/2^{N-1}) = 1/2 + F_{2v}^+(1/2)/4.
\end{array}
\end{equation}
In this case, the global maximum (respectively minimum) $M_v$ of the function $T_v$ on the segment $[0;1]$ 
can be calculated using the formula
\begin{equation} \label{eqMv_MvN}
M_v = \frac{M_{v, N}}{1-v^{N+1}}.
\end{equation}

2c) In the case of $v^{N+1}<0$ (that is, $v<0$ and $N$ even), the set $E_v$ consists of all points $x$ of the form
\begin{equation} \label{eqxExtrSoglPol2}
x = 1/2 \pm (P_{2v,N}(1/2)/4 + 1/2^{N+2} - y/2^{N+1}),
\end{equation}
where $y$ is any point of the global minimum (respectively maximum) of the function $T_v$ on the segment $[0;1]$.

In this case, the global maximum (respectively minimum) $M_v$ of the function $T_v$ on the segment $[0;1]$ 
can be calculated using the formula
\begin{equation} \label{eqMaxT_v5}
M_v = M_{v, N} + v^{N+1} m_v,
\end{equation}
where 
$m_v = \min_{y\in[0,1]}T_v(y)$ (respectively, $m_v = \max_{y\in[0,1]}T_v(y)$).
\end{Teor}
{\bf Proof} will be done for maxima only, since the proof for minima is similar.

1). Let's consider the first case, where the point $2v$ is consistent with a series (not polynomial)  
$F_{2v}(x) = c_0 + c_1 x + \ldots + c_n x^n + \ldots$.

1a). First, find the set $E_v\cap[0; 1/2]$. 
To do this, we build a sequence of nested segments $[a_n,b_n]$ containing it.
For any $n=0,1,\ld$ put $P_{n}(t) = c_0 + c_1 t + \ldots + c_n t^n$ and set the numbers $a_n = 1/2 - P_n(1/2)/4-1/2^{n+2}$, $b_n = a_n + 1/2^{n+1}$.
Hence, the next properties (1.0$_n$)--(1.2$_n$) follow from the inequality~\eqref{eqSoglKoef}:

(1.0$_n$) $b_n - a_n = 1/2^{n+1}$;

(1.1$_n$) $[a_n, b_n]$ has the form $[m/2^n-1/2^{n+1}, m / 2^n]$ or $[m/2^n, m / 2^n+1/2^{n+1}]$ for some $m\in\ZZ$;

(1.2$_n$) $c_{n+1}\cdot P_n(2v) < 0$.

Now let's prove by induction, that for any $n=0,1,\ld$ the following properties (1.3$_n$)--(1.5$_n$) also hold:

(1.3$_n$) $S'_{v, n}(x) = P_n(2v)$ for any $x$ in the interval $(a_n, b_n)$;

(1.4$_n$) $\bigl(E_{v, n}\cap[0; 1/2]\bigr) \subset [a_n, b_n]$ and $\bigl(E_v\cap[0; 1/2]\bigr) \subset [a_n, b_n]$;

(1.5$_n$) if $n\ges1$ then $[a_n, b_n] \subset [a_{n-1}, b_{n-1}]$.

In the case of $n=0$, we have: 
$a_n=0$, $b_n=1/2$, $S_{v, n}(x) = T_0(x) = x$ when $x\in[a_n, b_n]$, 
$a_{n+1} = (1-c_1)/8=1/4$, $b_{n+1} = (3-c_1)/8=1/2$.
Therefore, the properties (1.3$_n$)--(1.5$_n$) are met.

Now assume that the properties (1.3$_n$)--(1.5$_n$) are met, and prove that the properties (1.3$_{n+1}$)--(1.5$_{n+1}$) are met.
We have:
$a_{n+1} = 1/2 - P_{n+1}(1/2)/4 - 1/2^{n+3} = a_n + (1 - c_{n+1})/2^{n+3}$, $b_{n+1} = a_{n+1} + 1/2^{n+2} = b_n - (1 + c_{n+1})/2^{n+3}$.
From the properties (1.2$_n$) and (1.3$_n$), it follows that for any $x\in(a_n, b_n)$ 
the equality $\sign(S'_{v, n}(x)) = \sign(P_n(2v)) = - c_{n+1}$  is true.
Let's first consider the case when $c_{n+1} = -1$.
In this case, $a_{n+1} = a_n + 1/2^{n+2}$ and $b_{n+1} = b_n$.
Therefore, $[a_{n+1}, b_{n+1}]$ is the right half of the segment $[a_n, b_n]$, that means the property (1.5$_{n+1}$) is fulfilled.
In addition, for any $x$ in the interval $(a_{n+1}, b_{n+1})$, the equalities $T'_0(2^{n+1}x) = -1 = c_{n+1}$ and
$$
S'_{v, n+1}(x) = S'_{v, n}(x) + (2v)^{n+1} T'_0(2^{n+1}x) = P_n(2v) + (2v)^{n+1} c_{n+1} = P_{n+1}(2v)
$$
are true.
Therefore, the property (1.3$_{n+1}$) is also met.
Since $\sign(S'_{v,n}(x)) = - c_{n+1} >0$ in this case, then, by virtue of property (1.1$_n$) and item~b) of the lemma~\ref{LemSvnLin}, 
the property (1.4$_{n+1}$) is also fulfilled.
Thus, in the case $c_{n+1} = -1$, the properties (1.3$_{n+1}$)--(1.5$_{n+1}$) are proved.
In the case of $c_{n+1} = 1$, they can be proved in a similar way.

So now, the sequence $[a_n, b_n]$, $n=0,1,\ld$ with the properties (1.0$_n$) - (1.5$_n$) is constructed.
Then it follows from Cantor's theorem on nested segments, 
that the set $E_v\cap[0; 1/2]$ consists of a single point $x^-(v)$, and
$$
x^-(v) = \lim_{n\to\infty} a_n = \lim_{n\to\infty} (1/2 - P_n(1/2)/4-1/2^{n+2}) =
1/2 - F_{2v}(1/2)/4.
$$
Thus, the formula~\eqref{eqx^-(v)} is proved.
Then the formula~\eqref{eqx_n^-} follows from the next chain of equalities:
$$
\sum_{n=1}^{\infty} \frac{x_n^-}{2^n} = x^-(v) =
\frac{1}{2} - \frac{F_{2v}(1/2)}{4} = 
\sum_{n=1}^{\infty} \frac{1}{2^{n+1}} - \sum_{n=0}^{\infty} \frac{c_n}{2^{n+2}} =
\sum_{n=1}^{\infty} \frac{1-c_{n-1}}{2^{n+1}}.
$$

Now find the set $E_v\cap[1/2;1]$.
Due to the identity $T_v(x) = T_v(1-x)$, this set also consists of a single point $x^+(v) = 1 - x^-(v)$.
Therefore the formula~\eqref{eqx^+(v)} for $x^+(v)$ and the formula~\eqref{eqx_n^+} for its binary digits are true.

1b). Next, we prove the formulae for the global maxima $M_v$ of the functions $T_v$.
Since $T_0$ has the period $1$, $T_0(x)=x$ for $x\in[0; 1/2]$ and $T_0(x)=1-x$ for $x\in[1/2;1]$, 
then for any $n=0,1,\ld$ the following equality are met:
$$
T_0(2^n x^-(v)) = T_0(x_1^-\ld x_n^-{.}x_{n+1}^- x_{n+2}^- \ld) = T_0(0{.}x_{n+1}^-x_{n+2}^-\ld) =
$$
$$
= \left\{\begin{array}{rl}
  0{.}x_{n+1}^-x_{n+2}^-\ld&\text{ when }x_{n+1}^- =0\\
1-0{.}x_{n+1}^-x_{n+2}^-\ld&\text{ when }x_{n+1}^- =1\\
\end{array}\right. =
$$
\begin{equation} \label{eqT0(2nx-(v))}
= x_{n+1}^- + (1-2x_{n+1}^-)\cdot 0{.}x_{n+1}^-x_{n+2}^-\ld =
x_{n+1}^- + 2^n(1-2x_{n+1}^-)\sum_{k=n+1}^\infty \frac{x_k^-}{2^k}.
\end{equation}
By virtue of the formula~\eqref{eqx_n^-}, $x_n^- =(1-c_{n-1})/2$ for any $n\in\NN$, so
$$
T_0(2^n x^-(v)) = \frac{1-c_{n}}{2} + 2^n c_n\sum_{k=n+1}^\infty \frac{1-c_{k-1}}{2^{k+1}} =
\frac12-2^n c_n\sum_{i=n}^\infty \frac{c_i}{2^{i+2}}.
$$
From here we get the proved formula~\eqref{eqMaxT_v1}:
$$
T_v(x^-(v)) = \sum_{n=0}^\infty v^n T_0(2^n x^-(v)) = 
\sum_{n=0}^\infty v^n\Bigl(\frac12-2^n c_n\sum_{i=n}^\infty \frac{c_i}{2^{i+2}}\Bigr) =
$$
$$
= \frac{1}{2(1-v)} - \frac14\sum_{n=0}^\infty c_n(2v)^n \sum_{i=n}^\infty \frac{c_i}{2^i}.
$$
Now let's proof the formula~\eqref{eqMaxT_v2}.
Take any number $r$ from the interval $(\max (1/2, v); 1)$.
Then, applying the following formula for the coefficients of the Taylor series (see~\cite[Ch.~V, \S2, p.~1, p.~196]{Privalov}):
\begin{equation} \label{eqFormKowi}
c_n = \frac{F_{2v}^{(n)}(0)}{n!} = \frac{1}{2\pi i}\int_{|z|=r}\frac{F_{2v}(z)}{z^{n+1}} dz,
\end{equation}
we find:
$$
\sum_{i=n}^\infty \frac{c_i}{2^i} = \frac{1}{2\pi i}\int_{|z|=r}\sum_{i=n}^\infty \frac{F_{2v}(z)}{2^i z^{i+1}} dz =
\frac{1}{\pi i}\int_{|z|=r} \frac{F_{2v}(z)}{(2z)^n(2z-1)} dz.
$$
Substituting this expression to~\eqref{eqMaxT_v1}, we get the equality~\eqref{eqMaxT_v2}:
$$
T_v(x^-(v)) = \frac{1}{2(1-v)} - 1/4\sum_{n=0}^\infty c_n(2v)^n \frac{1} {\pi i}\int_{|z|=r} \frac{F_{2v}(z)}{(2z)^n(2z-1)} dz =
$$
$$
= \frac{1}{2(1-v)} - \frac{1}{4\pi i}\int_{|z|=r} \frac{F_{2v}(z)}{2z-1} \sum_{n=0}^\infty c_n(v/z)^n dz =
$$
$$
= \frac{1}{2(1-v)} - \frac{1}{4\pi i}\int_{|z|=r} \frac{F_{2v}(z)F_{2v}(v/z)}{2z-1} dz.
$$

2). Let's move to the second case, where the point $2v$ is consistent 
with the polynomial (not series) $P_{2v,N}(x) = c_0 + c_1 x + \ldots + c_N x^N$.

2a) As in the case 1) (when $2v$ was consistent with a series), we construct a sequence (at this time, finite) 
of nested segments $[a_n, b_n]$ with properties (1.0)$_n$-(1.5)$_n$ for any $n=0,1,\ld, N$.
In particular, the equalities~\eqref{eqabN} are true: $a_N = 1/2 - P_{2v,N}(1/2)/4 - 1/2^{N+2}$, $b_N = a_N + 1/2^{N+1}$.

Since $P_{2v, N}(2v) = 0$, then, by virtue of the property (1.3$_N$), $S'_{v, N} = 0$ on the interval $(a_N, b_N)$.
Therefore, the function $S_{v, N}$ is constant on the segment $[a_N, b_N]$.
By property (1.4$_N$), the sets $E_{v, N}\cap[0; 1/2]$ and $E_v\cap[0;1/2]$ lie on the segment $[a_N,b_N]$.
So $E_{v, N}\cap[0; 1/2] = [a_N, b_N]$.
Hence, due to the symmetry of the function $S_{v, N}$, we get the equality $E_{v, N} = [a_N, b_N]\cup[1-b_N, 1-a_N]$.

From the property (1.1$_N$), it follows that $a_N$ has a finite binary expansion $a_N = 0{.}x_1^-x_2^- \ld x_{N+1}^-$, 
the points of the segment $[a_N,b_N]$ have a binary expansion of the form $y = 0{.}x_1^-x_2^- \ld x_{N+1}^-\ld$, 
and the points of $[1-b_N,1-a_N]$ have binary expansion $z = 0{.}x_1^+x_2^+ \ld x_{N+1}^+\ld$, 
where $x_i^+ = 1-x_i^-$ for $i=1,\ld,N+1$, and $1$ in the period is allowed in binary expansions of the points $y$ and $z$.
Thus, the set $E_{v, N}$ consists of all points $x$ whose binary expansion has the form
\begin{equation} \label{eqxEvN}
x = 0{.}\left[\begin{array}{l}
x_1^+x_2^+\ld x_{N+1}^+\\x_1^-x_2^-\ld x_{N+1}^-
\end{array}\right. \ld.
\end{equation}
The relations~\eqref{eqx_n^-Pol} follow from the formulae $x_i^+ = 1-x_i^-$ for $i=1,\ld,N+1$ and the equalities
$$
\sum_{n=1}^{N+1} \frac{x_n^-}{2^n} = a_N = \frac12 - \frac{1}{2^{N+2}} - \frac14 P_{2v, N}\bigl(\frac12\bigr) =
\sum_{n=1}^{N+1} \frac{1-c_{n-1}}{2^{n+1}}.
$$

Let's go to the proof of the formulae \eqref{eqMvN1} and \eqref{eqMvN2}.
Similar to the formula~\eqref{eqT0(2nx-(v))}, for any $n=0,1,\ld, N$ we get the equality
$$
T_0(2^n a_N) = (1-2x_{n+1}^-)2^n\sum_{k=n+1}^{N+1} x_k^-/2^k + x_{n+1}^-.
$$
Therefore, since $x_n^- =(1-c_{n-1})/2$ by virtue of the formula~\eqref{eqx_n^-Pol}, we find:
$$
T_0(2^n a_N) = 1/2 - c_n 2^n\sum_{i=n}^N c_i/2^{i+2} - c_n2^n/2^{N+2}.
$$
Hence, in view the equality $P_{2v, N}(2v)=0$, the formula~\eqref{eqMvN1} follows:
$$
M_{v, N} = T_v(a_N) = \sum_{n=0}^\infty v^n T_0(2^n a_N) = \sum_{n=0}^N v^n T_0(2^n a_N) =
$$
$$
= \sum_{n=0}^N v^n/2 - \sum_{n=0}^N c_n(2v)^n\sum_{i=n}^N c_i/2^{i+2} - P_{2v, N}(2v)/2^{N+2} =
$$
$$
= (1-v^{N+1})/(2(1-v)) - 1/4\sum_{n=0}^N c_n(2v)^n \sum_{i=n}^N c_i/2^i.
$$

Now we prove the formula~\eqref{eqMvN2}. Let $r\in(0; 1/2)\cup(1/2;1)$.
Applying the formula~\eqref{eqMaxT_v2}, as in the proof of equality~\eqref{eqFormKowi}, we get:
$$
\sum_{i=n}^N c_i/2^i = \frac{1}{2\pi i}\int_{|z|=r}\sum_{i=n}^N \frac{F_{2v}(z)}{2^i z^{i+1}} dz =
\frac{1}{\pi i}\int_{|z|=r} \frac{F_{2v}(z)}{2z-1}\Bigl(\frac{1}{(2z)^n} - \frac{1}{(2z)^{N+1}}\Bigr) dz.
$$
Substituting this expression into~\eqref{eqMvN1}, swapping the sum and integral and 
taking into account the equality $F(2v)=0$, we come to the formula~\eqref{eqMvN2}:
$$
M_{v, N} = \frac{1-v^{N+1}}{2(1-v)} - \frac{1}{4}\sum_{n=0}^N c_n(2v)^n 
\frac{1}{\pi i}\int_{|z|=r} \frac{F_{2v}(z)}{2z-1}\Bigl(\frac{1}{(2z)^n} - \frac{1}{(2z)^{N+1}}\Bigr) dz =
$$
$$
= \frac{1-v^{N+1}}{2(1-v)} - \frac{1}{4\pi i}
\int_{|z|=r} \frac{F_{2v}(z)}{2z-1}\Bigl(\sum_{n=0}^N c_n\Bigl(\frac{v}{z}\Bigr)^n - 
\frac{F(2v)}{(2z)^{N+1}}\Bigr) dz =
$$
$$
= \frac{1-v^{N+1}}{2(1-v)} - \frac{1}{4\pi i}
\int_{|z|=r} \frac{F_{2v}(z)}{2z-1}F(v/z) dz.
$$

2b) Let $v^{N+1}>0$, that is, $v>0$ or $N$ is odd.
Write $T_v(x)$ like this:
$$
T_v(x) = \sum\limits_{k=0}^\infty \sum\limits_{n=k(N+1)}^{k(N+1)+N} v^n T_0(2^nx) =
\sum\limits_{k=0}^\infty \sum\limits_{i=0}^{N} v^{k(N+1)+i} T_0(2^{k(N+1)+i}x) =
$$
$$
= \sum\limits_{k=0}^\infty v^{k(N+1)} \sum\limits_{i=0}^{N} v^i T_0(2^i 2^{k(N+1)}x) =
\sum\limits_{k=0}^\infty v^{k(N+1)} S_{v, N}(2^{k(N+1)}x).
$$
Hence the estimation follows:
$T_v(x) \les 
\sum\limits_{k=0}^\infty v^{k(N+1)} M_{v, N} = M_{v, N}/(1-v^{N+1})$.
Equality is achieved in this estimation only if ${S_{v, N}(2^{k(N+1)}x) = M_{v, N}}$ for any $k=0,1,2,\ld$.
Let $x = 0{.}x_1 x_2 \ld$. Then, due to the periodicity of the function $S_{v, N}$, we have:
$$
S_{v, N}(2^{k(N+1)}x) = S_{v, N}(x_1 x_2 \ld x_{k(N+1)}{.}x_{k(N+1)+1} x_{k(N+1)+2} \ld) =
$$
$$
= S_{v, N}(0{.}x_{k(N+1)+1} x_{k(N+1)+2} \ld).
$$
So $S_{v, N}(2^{k(N+1)}x) = M_{v, N}$ only if $0{.}x_{k(N+1)+1} x_{k(N+1)+2}\ld \in E_{v, N}$.
The last condition is equivalent by the formula~\eqref{eqxEvN} to the fact that 
the set $\{x_{k(N+1)+1}, x_{k(N+1)+2},\ld, x_{k(N+1)+N+1)}\}$ coincides with one of two following sets: 
$\{x_1^-, x_2^-,\ld, x_{N+1)}^-\}$ or $\{x_1^+, x_2^+,\ld, x_{n+1)}^+\}$.
Therefore, the function $T_v$ reaches a global maximum on the segment $[0;1]$ at the point $x$ if and only if 
its binary expansion has the form~\eqref{eqxExtrSoglPol1}.
In this case, $M_v = M_{v, N}/(1-v^{N+1})$, that is, we proved the equality~\eqref{eqMv_MvN}.

It follows from the proven formula~\eqref{eqxExtrSoglPol1} that
$$
\inf E_v = 0{.}x_1^-x_2^- \ld x_{N+1}^- \ x_1^-x_2^- \ld x_{N+1}^- \ \ld =
a_N/(1-1/2^{N+1}).
$$
Therefore, taking into account the formulae \eqref{eqSumF+-} and~\eqref{eqabN}, we get the equality~\eqref{eqinfEv}.
Next, according to the item~b) of the lemma~\ref{LemSvnLin}, 
the set $E_v\cap[0; 1/2]$ is symmetric about the middle of $[a_N, b_N]$, 
so $\sup(E_v\cap[0; 1/2]) = a_N+b_N - \inf(E_v\cap[0; 1/2])$.
From here, and from the equalities \eqref{eqSumF+-}, \eqref{eqabN}, \eqref{eqinfEv}, we get the formula~\eqref{eqinfEv012}.
The formula~\eqref{eqsupEv} follows from~\eqref{eqinfEv} and from the equality $\sup E_v = 1 - \inf E_v$.

2c) Let $v^{N+1}<0$, that is, $v<0$ and $N$ even.
Then, using the functional equation~\eqref{eqFunkUr}, we rewrite $T_v(x)$ in the following form:
\begin{equation} \label{eqTv_v<0}
T_v(x) = S_{v, N}(x) + v^{N+1}T_v(2^{N+1}x).
\end{equation}
By virtue of item~2a) we have: $E_v\cap[0; 1/2] \subset E_{v, N}\cap[0; 1/2] = [a_N, b_N]$.
Therefore, the set $E_v\cap[0; 1/2]$ consists of those points $x$ of the segment $[a_N, b_N]$, 
in which the function $T_v(2^{N+1}x)$ is minimal.
These points have the form $x = a_N + y/2^{N+1}$, where $y$ is the point of the lowest value of the function $T_v$ on $[0;1]$.
From here, taking into account the formula~\eqref{eqabN} and the symmetry of the set $E_v$ about the point $1/2$, 
we get that $E_v$ consists of all points of the form~\eqref{eqxExtrSoglPol2}.

From the constancy of $S_{v, N}$ on the segment $[a_N, b_N]$ and the equality~\eqref{eqTv_v<0}, it follows that
$$
M_v = M_{v, N} + v^{N+1}\min_{x\in[a_N, b_N]}T_v(2^{N+1}x) =
M_{v, N} + v^{N+1}\min_{y\in[2^{N+1}a_N, 2^{N+1}b_N]}T_v(y).
$$
Hence, given that the length of the segment $[2^{N+1}a_N,2^{N+1}b_N]$ is equal 
to the period of the function $T_v$, i.e. $1$, we get the equality~\eqref{eqMaxT_v5}. 
The theorem is proved. 
\qed 

\medskip

Note that although the maximum of the $T_v$ function is expressed through its minimum and vice versa
in some cases in the theorem~\ref{TeorPoiskGlobExtr}, there is no vicious circle here. 
This will be seen from the theorem~\ref{TeorGlobMin} of the next section.

\begin{Prim} \label{PrimMaxV1/2sqrt2}
For $v = 1/(2\sqrt{2}) \approx 0.353553$, the $T_v$ function has two maximum points on the segment $[0,1]$: $5/12 \approx 0{.}41667$ and $7/12 \approx 0{.}58333$.
In this case, $M_v = (26+3\sqrt2)/56$.
\end{Prim}
{\bf Proof.}
As shown in the example~\ref{PrimSoglW1/sqrt2}, the point $2v = 1/\sqrt{2}$ is consistent with the series $F_{2v}$,
that have the sum $F_{2v}(x) = (1-2x^2)/(1+x)$ for $|x|<1$. 
Then, using the formulae \eqref{eqx^-(v)} and~\eqref{eqx^+(v)}, 
we get the maximum points of $T_v$ on the segment $[0, 1]$: 
$x_{max}^- = 1/2-1/4\cdot F_{2v}(1/2) = 5/12$ and $x_{max}^+ = 1-x_{max}^- = 7/12$.
Next, using the formula~\eqref{eqMaxT_v2} with $r=3/4$ we have:
$$
M_v = \frac{1}{2(1-v)} - \frac{1}{4\pi i}\int_{|\xi|=3/4} 
\frac{F_{2v}(\xi)F_{2v}(v/\xi)}{2\xi-1}d\xi =
$$
$$
= \frac{1}{2(1-v)} - \frac{1}{4\pi i}\int_{|\xi|=3/4} 
\frac {(1-2\xi^2)(1-2v^2/\xi^2)}{(2\xi-1)(1+\xi)(1+v/\xi)}d\xi =
$$
$$
= \frac{2}{4-\sqrt2} + \frac{1}{4\pi i}\int_{|\xi|=3/4}
\frac{(2\xi^2-1)(2\xi+1)}{\xi(\xi+1)(4\xi+\sqrt2)} d\xi.
$$
Calculating the integral, we get the equality $M_v = (26+3\sqrt2)/56$.
\qed


\section{\boldmath Global minima of the functions $T_v$ in the general case $v\in(-1;1)$}
\label{SectGlobMin}

In the single theorem~\ref{TeorGlobMin} of this section, 
we find the global minima of the functions $T_v$ for all points $v\in(-1;1)$, 
using the theorem~\ref{TeorPoiskGlobExtr}.
\begin{Teor}
\label{TeorGlobMin}
The following statements are true:

1) If $v\in(-1/2;1)$, then the global minimum of the function $T_v$ on the segment $[0;1]$ is equal to $0$
and is achieved only at two points: $x_{min}^- = 0$ and $x_{min}^+ = 1$.

2) (Cf. with section~6 of ~\ref{SectIzvSvoistva}) If $v = -1/2$, then the set of points of the global minimum of the function $T_v$ on the segment $[0;1]$ has a Hausdorff dimension $1/2$ and consists of all points $x_{min}$ with a binary expansion of the form
\begin{equation} \label{eqxmin(-1/2)}
x_{min} = 0{.}
\left[\begin{array}{l}00\\11\end{array}\right.
\left[\begin{array}{l}00\\11\end{array}\right. \ld,
\end{equation}
that is, having the expansion $x_{min} = 0{.}y_1y_2\ldots$, where $y_1, y_2,\ldots\in\{0;3\}$, in a system with base 4.
In this case, $\min_{x\in[0;1]}T_v(x) = 0$.

3) If $v\in(-1; -1 / 2)$, then the function $T_v$ reaches the global minimum on $[0;1]$ at only two points: 
$x_{min}^- = 1/5$ and $x_{min}^+ = 4/5$.
In this case, $\min_{x\in[0;1]}T_v(x) = T_v(1/5) = (1+2v)/(5(1-v^2))$.
\end{Teor}
{\bf Proof.} 
1). If $v\in(-1/2; 1)$, then $2v\in(-1;2)$. 
Therefore, it follows from item~1) of the theorem~\ref{TeorVseAntisogl} that the point $2v$ is anti-consistent with the series $A (x)=1+x+x^2+\ldots$, and $A (x)=1/(1-x)$ for $|x|<1$.
Due to item~1a) of the theorem~\ref{TeorPoiskGlobExtr}, there are two minimum points, which can be calculated using the formulae \eqref{eqx^-(v)} and \eqref{eqx^+(v)}: $x_{min}^- = 1/2 - A(1/2)/4 = 0$, $x_{min}^+ = 1 - x_{min}^- = 1$.

2) If $v = -1/2$, then $2v = -1$.
By the theorem~\ref{TeorVseAntisogl}, the point $2v$ is anti-consistent with the polynomial $P(x)=1+x$.
Then according to item~2) of the theorem~\ref{TeorPoiskGlobExtr} we have:

2a) Since $N=1$, then $a_N(v) = 1/2-P(1/2)/4-1/2^{N+2} = 0$.
In this case, the binary expansion of the length $N+1$ of the point $a_N(v)$ has the form $a_N(v) = 0{.}00$.

2b) Since $N=1$ is odd, then the set $E_v$ has Hausdorff dimension $1/(N+1)=1/2$ and consists of all points $x_{min}$,
having the binary expansion of the form~\eqref{eqxExtrSoglPol1}, where $x_1^-=x_2^- =0$, $x_1^+=x_2^+ = 1-x_1^- =1$. 
Therefore we get the formula~\eqref{eqxmin(-1/2)}.
Since $0\in E_v$, then $\min_{x\in[0;1]}T_v(x) = T_v(0) = 0$.

3) If $v\in(-1; -1/2)$, then $2v\in(-2;-1)$.
So, according to item~3) of the theorem~\ref{TeorVseAntisogl}, the point $2v$ is anti-consistent with the series $A (x)=1+x-x^2-x^3+x^4+x^5-\ldots$, and $A (x)=(1+x)/(1+x^2)$ for $|x|<1$.
According to item~1a) of the theorem~\ref{TeorPoiskGlobExtr}, 
there are two points of the global minimum of the function $T_v$ on the segment $[0;1]$.
These points can be calculated using the formulae \eqref{eqx^-(v)} and~\eqref{eqx^+(v)}: 
$x_{min}^- = 1/2 - A(1/2)/4 = 1/5$ and $x_{min}^+ = 1 - x_{min}^- = 4/5$.

To calculate the minimum value of $T_v (1/5)$, we could use the formula~\eqref{eqMaxT_v2}, 
but instead we will use the functional equation~\eqref{eqFunkUr} for $x=1/5$ and $N=1$: 
$T_v (1/5) = T_0(1/5) + v T_0(2/5) + v^2 T_v(4/5)$.
Hence, by virtue of the equalities $T_0 (1/5) =1/5$, $T_0(2/5) = 2/5$, and $T_v(4/5) = T_v(1/5)$, 
we get the equation $T_v (1/5) = 1/5 + 2v/5 + v^2 T_v(1/5)$. Solving it with respect to $T_v(1/5)$, 
we get the required formula $T_v(1/5) = (1+2v)/(5(1-v^2))$.
\qed 

\section{\boldmath Global maxima of the functions $T_v$ in the general case $v\in(-1;1)$}
\label{SectGlobMaks}

In this section, we study the global maxima $M_v$ of the functions $T_v$ on the segment $[0;1]$, 
as well as the set $E_v$ of points of the global maximum, for $v\in(-1;1)$.
In the case where $v\in(-1; 1/4]\cup[1/2;1)$, the situation is much simpler than in the case of $v\in (1/4;1/2)$.
The first case is fully discussed in the first subsection of this section.
All further subsections of this section are devoted to the second case.

\subsection{\boldmath Global maxima of the functions $T_v$ in the case when $v\in(-1; 1/4]\cup[1/2;1)$}

Here, as in the section~\ref{SectGlobMin}, we will base on the theorem~\ref{TeorPoiskGlobExtr}.
\begin{Teor}
\label{TeorPoiskGlobMax}
The following statements take place:

1) If $1/2 < v < 1$, then $E_v = \{1/3, 2/3\}$ and $M_v = 1/(3(1-v))$.

2) If $v = 1/2$, then the set $E_v$ consists of points $x_{max}$ having a binary expansion $x_{max} = 0, x_1 x_2\ldots x_n\ldots$ satisfying the condition $x_{2k+1}+x_{2k+2}=1$ for $k=0,1,\ldots$ (see\ \cite[Lieu~1]{Kahane} or the theorem~\ref{TeorKahane}).
These points also can be written as
\begin{equation} 
x_{max} = 0{.}
\left[\begin{array}{l}01\\10\end{array}\right.
\left[\begin{array}{l}01\\10\end{array}\right. \ld.
\end{equation}
Also, $M_v = 2/3$ (see~\cite[Prop.~4.2]{Kruppel2007}), 
and the next equalities are true: $\inf E_v = 1/3$, $\sup(E_v\cap[0;1/2]) = 5/12$, $\sup E_v = 2/3$.

3) If $-1/2 \leq v \leq 1/4$, then $E_v=\{1/2\}$ and $M_v = 1/2$.

4) If $-1 < v < -1/2$, then

4a) for any $k\in {\mathbb N}$ polynomial $\widetilde{P}_{2k}(x) = 1-2x-4x^2 - \ld-2^{2k}x^{2k}$ has a single negative root $x=v_k$, and $v_k$ belongs to the interval $(-1;-1/2)$. 
The sequence of roots $\{v_k\}_{k=1}^\infty$ strictly increases, 
and the asymptotic equality $v_k = -1/2-(\ln3)/(4k)+\underline{O}(1/k^2)$ is satisfied for $k\to\infty$.

4b) If $k\in {\mathbb N}$ and $v\in(v_{k-1}, v_k)$ (here $v_0=-1$), then the set $E_v$ consists of two points: 
$E_v = \{x_{max, k}^-, x_{max, k}^+\}$, where
\begin{equation} \label{eqxmax_kpm}
x_{max, k}^- = \frac{1}{2} - \frac{1}{5\cdot2^{2k-1}}; \quad
x_{max, k}^+ = \frac{1}{2} + \frac{1}{5\cdot2^{2k-1}}.
\end{equation}
Besides,
\begin{equation} \label{eqmaxTvk}
M_v = \frac{1}{2} + \frac{4v-1}{5\cdot2^{2k-1}(1-2v)} - \frac{3v^{2k+1}}{5(1-v^2)(1-2v)}.
\end{equation}

4c) If $k\in{\mathbb N}$, and $v = v_k$, then $E_v = \{x_{max,k}^-, x_{max,k+1}^-, x_{max,k+1}^+, x_{max,k}^+\}$, 
where the points $x_{max,k}^\pm$, $x_{max,k+1}^\pm$ are defined by the formulae~\eqref{eqxmax_kpm}.
The value of $M_v$ can be calculated using the formula~\eqref{eqmaxTvk}.
\end{Teor}
{\bf Proof.}
1). If $1/2<v<1$, then $1<2v<2$. 
So, according to item~1) of the theorem~\ref{TeorSoglW<1/2>1}, the point $2v$ is consistent with the series, 
that has the sum $F(x)=1/(1+x)$ for $|x|<1$.
According to item~1a) of the theorem~\ref{TeorPoiskGlobExtr}, 
the set $E_v$ consists of two points $x_{max}^- $ and $x_{max}^+$, 
which can be calculated using the formulae \eqref{eqx^-(v)} and \eqref{eqx^+(v)}: 
$x_{max}^- = 1/2-F(1/2)/4 = 1/3$, $x_{max}^+ = 1 - x_{max}^- = 2/3$.

To calculate the value of $M_v = T_v(1/3)$, we apply the functional equation~\eqref{eqFunkUr} for $x=1/3$ and $N=0$: 
$T_v(1/3) = T_0 (1/3) + v T_v(2/3)$.
Hence, by virtue of the equalities $T_0 (1/3) =1/3$ and $T_v(2/3) = T_v(1/3)$, 
we get the equation $T_v (1/3) = 1/3 + v T_v(1/3)$.
So we find from it: $M_v = T_v(1/3) = 1/(3(1-v))$.

2) For $v = 1/2$, the structure of the set $E_v$ is described by Kahane (see theorem~\ref{TeorKahane} or~\cite[Lieu~1]{Kahane}).
The value $M_v$ is calculated in~\cite[Prop.~4.2]{Kruppel2007}.
To search for $\sup(E_v\cap[0; 1/2])$, note that the point $2v=1$ is consistent with the polynomial $P(x)=1-x$ of degree $N=1$.
Therefore, according to the formulae \eqref{eqinfEv}, \eqref{eqinfEv012} and\eqref{eqsupEv} of the theorem~\ref{TeorPoiskGlobExtr} we have: $\inf E_v = 1/2 - P(1/2)/(4-1) = 1/3$, $\sup(E_v\cap[0;1/2]) = 1/2 - 1/2\cdot P(1/2)/3 = 5/12$, $\sup E_v = 1 - \inf E_v = 2/3$.

3) If $-1/2 \leq v \leq 1/4$, then $-1 \leq 2v \leq 1/2$.
So, according to item~3) of the theorem~\ref{TeorSoglW<1/2>1}, the point $2v$ is consistent with a series, 
which has the sum $F(x) = (1-2x)/(1-x)$ for $|x|<1$.
According to item~1a) of the theorem~\ref{TeorPoiskGlobExtr}, there are two maximum points
(they coincide in this case): $x_{max}^- = 1/2 - F(1/2)/4 = 1/2$ and $x_{max}^+ = 1 - x_{max}^- = 1/2$.
Moreover, obviously, $M_v = T_v(1/2)=1/2$.

4) If $-1 < v < -1/2$, then $-2 < 2v < -1$.

4a) The statement of this item follows from item~4a) of the theorem~\ref{TeorSoglW<1/2>1}, 
where it is necessary to set $t=2x$. Then we get $v_k = w_k/2$ for all $k\in\NN$.

4b) According to item~4b) of the theorem~\ref{TeorSoglW<1/2>1}, for any $k\in\NN$ and $v\in(v_{k-1},v_k)$ (where $v_0=-1$), 
the point $2v$ is consistent with the series $F_k$, having the sum $F_k(x) = (1-2x)(1-x) + 2x^{2k+1}/((1-x)(1+x^2))$
for $|x|<1$.
According to item~1a) of the theorem~\ref{TeorPoiskGlobExtr}, 
$E_v$ consists of two maximum points $x_{max}^- $ and $x_{max}^+$, 
which can be calculated by the formulae \eqref{eqx^-(v)} and \eqref{eqx^+(v)}: 
$x_{max, k}^- = 1/2 - F_k(1/2)/4 = 1/2 - 1/(5\cdot2^{2k-1})$, 
$x_{max, k}^+ = 1 - x_{max}^- = 1/2 + 1/(5\cdot2^{2k-1})$.

The value $M_v$ can be calculated using the formula~\eqref{eqMaxT_v2} of the theorem~\ref{TeorPoiskGlobExtr}:
$$
M_v = \frac{1}{2(1-v)} - \frac{1}{4\pi i}\int_{|z|=r} \frac{F_k(z)F_k(v/z)}{2z-1}dz =
$$
$$
= \frac{1}{2(1-v)} - \frac{1}{4\pi i}\int_{|z|=r} 
\Bigl( \frac1{z-1} + \frac{2z^{2k+1}}{(2z-1)(1-z)(1+z^2)} \Bigr) \cdot
$$
$$
\cdot \Bigl( \frac{z-2v}{z-v} + \frac{2v^{2k+1}}{z^{2k-2}(z-v)(z^2+v^2)} \Bigr) dz.
$$
Calculating this integral using residues, we get the necessary formula~\eqref{eqmaxTvk}.

4c) If $v = v_k$ for some $k\in {\mathbb N}$, then, according to item~4b) of the theorem~\ref{TeorSoglW<1/2>1}, the point $2v$ is consistent with the polynomial $P_{2k}(x) = 1-x-x^2 - \ld-x^{2k}$ of degree $N=2k$.
For $x\neq1$, we have: $P_{2k}(x) = (1-2x+x^{2k+1})/(1-x)$.
So, according to item~2b) of the theorem~\ref{TeorPoiskGlobExtr}, 
the set $E_v$ consists of all points $x_{max}$ of the form~\eqref{eqxExtrSoglPol2}:
$$
x_{max} = 1/2 \pm (P_{2k}(1/2)/4 + 1/2^{2k+2} - y/2^{2k+1}) =
$$
$$ 
= 1/2 \pm \Bigl( (1-\sum_{n=1}^{2k}1/2^n)/4 + 1/2^{2k+2} - y/2^{2k+1}\Bigr) =
1/2 \pm (1-y)/2^{2k+1},
$$
where $y$ is any point of global minimum of the function $T_v$ on $[0;1]$.
According to item~3) of the theorem~\ref{TeorGlobMin}, 
the global minimum is reached only at~two points  in this case: $y^- = 1/5$ and $y^+ = 4/5$.
Therefore, $E_v$ consists of two pairs of points: 
the first pair is $1/2 \pm (1-1 / 5)/2^{2k+1} = 1/2 \pm 1/(5\cdot2^{2k-1})$ 
and the second pair is $1/2 \pm (1-4/5)/2^{2k+1} = 1/2 \pm 1/(5\cdot2^{2k+1})$.
According to the formulae~\eqref{eqxmax_kpm}, 
these points coincide with $x_{max, k}^\pm$ and $x_{max, k+1}^\pm$ respectively.
That's what we needed to prove.
\qed

\medskip
An example of the application of the proved theorem is given in the subsection~\ref{PrimerMaks}.


\subsection{On rational points of the global maximum}

\begin{Teor}
1) Suppose that $v\in(1/4;1/2)$, the point $2v$ is consistent with the series $F (x)=c_0+c_1x+c_2x^2+\ld$, 
and the point of the global maximum of the function $T_v$ on $[0; 1/2]$ is rational
(it is the only one by the item~1a of the theorem~\ref{TeorPoiskGlobExtr}). 
Then $F(x)$ can be represented as the following ratio of two polynomials for any $|x|<1$:
\begin{equation} \label{eqFrat}
F(x)=\frac{c_0+c_1x+\ld+c_{m-1}x^{m-1}+b_0x^m+b_1x^{m+1}+\ld+b_nx^{m+n}}{1-x^m},
\end{equation}
where $m\in\NN$; $n\in\{0,1,2,\ld\}$; $b_0, b_1,\ld,b_n\in\{-2,0,2\}$; $c_0, c_1,\ld,c_{m-1}\in\RR$.

2) Suppose that $v\in(-1;1)$ and some point of the global maximum of the function $T_v$ on $[0;1]$ is binary rational.
Then $v\in[-1/2;1/4]$, $E_v=\{1/2\}$ and $M_v=1/2$.
\end{Teor}
{\bf Proof.}
1) Let the global maximum point $x_{max}$ of the function $T_v$ on the segment $[0; 1/2]$ have a binary expansion $x_{max} = 0{.}x_1x_2\ld$.
Then $x_k = (1-c_{k-1})/2$ for all $k\in\NN$, according to the formula~\eqref{eqx_n^-}.
So due to the rationality of $x_{max}$, the sequence $c_0, c_1,c_2,\ld$ is periodic with some period $m$ starting from some number $n$.
Therefore, for any $|x|<1$ the function $F(x)$ can be represented as
$$
F(x) = c_0+c_1x+\ld+c_nx^n + \sum_{k=0}^\infty (c_{n+1}x^{n+km+1}+\ld+c_{n+m}x^{n+km+m}) =
$$
$$
= c_0+c_1x+\ld+c_nx^n + \frac{c_{n+1}x^{n+1}+\ld+c_{n+m}x^{n+m}}{1-x^m} =
$$
$$
= \frac{c_0+c_1x+\ld+c_{m-1}x^{m-1} + 
(c_{m}-c_0)x^{m}+\ld+(c_{n+m}-c_n)x^{n+m}}{1-x^m}.
$$
Taking $b_i=c_{i+m}-c_i\in\{-2,0,2\}$ for $i=0,1,\ld,n$, we get the equality~\eqref{eqFrat}.

2) If the point $x_{max} = 0{.}x_1x_2\ld$ is a binary rational point of the global maximum, 
then $x_n=x_{n+1}$ for any $n\ges N$ for some $N\in\NN$.
Due to the items 1), 2) and 4) of the theorem~\ref{TeorPoiskGlobMax}, this is not possible for $v\in(-1;-1/2)\cup[1/2;1)$.
Therefore $v\in[-1/2; 1/2)$.
We have $E_v=\{1/2\}$ and $M_v=1/2$ for any $v\in[-1/2;1/4]$, by virtue of item~3) of the theorem~\ref{TeorPoiskGlobMax}.
So it remains to show that the case of $v\in (1/4;1/2)$ is also impossible.
Let's consider two sub-cases.

2a) If the point $2v$ is consistent with the series $F(x)=c_0+c_1x+c_2x^2+\ld$, 
then $x_n = (1-c_{n-1})/2$ for any $n\in\NN$, according to the formula~\eqref{eqx_n^-}.
Therefore, $c_{n-1}=c_n$ for $n\ges N$.
So for any $m\in\NN$ the inequality $1-2v-\ld-(2v)^m \ges 0$ is satisfied, by virtue of the lemma~\ref{LemOdinZnak}.
If $m\to\infty$, then we get an estimation $1-2v-(2v)^2 - \ld = (1-4v)/(1-2v) \ges 0$.
Hence, $v\les1/4$.

2b) If $2v$ is consistent with the polynomial $P(x)=c_0+c_1x+\ld+c_Nx^N$, then the formula~\eqref{eqxExtrSoglPol1} implies that $x_1=\ld=x_{N+1}$.
So $c_0=\ld=c_N$ according to the formula~\eqref{eqx_n^-Pol}.
But this is also impossible, since $c_0=1$ and $c_1=-1$ by the item~1) of the remark~{ZamSoglFun}.
The theorem is proved.
\qed

\subsection{\boldmath How do the maximum points of the function $T_v$ change when replacing $v$ with $\sqrt{v/2}$ or $2v^2$}

Before answering the question contained in the title of this subsection, 
let's define and study the properties of the number $\chi$ and the mapping of $H$ on $[0;1)$.

Further if $x$ is a number, then $\overline{x}$ will denote the number $1-x$.
\begin{Opr} \label{OprCHi}
Let's define $\chi$ as a real number that has the binary expansion $\chi = 0{.}x_1 x_2 \ld$, 
where $x_1=0$, and the subsequent binary digits $x_2, x_3,\ld$ are calculated 
using the recurrent relations $x_{2k}=\overline{x_k}=1-x_k$ and $x_{2k+1}=x_{k+1}$ for each $k\in\NN$.
\end{Opr}

\begin{Zam} 
1) By definition~\ref{OprCHi}, the binary and decimal expansions of the number $\chi$ begin as follows: 
$\chi = 0{.}0110100110010110\ld_2 = 0{.}412454033\ld_{10}$.

2) The sequence of binary digits of the number $\chi$ 
satisfies the equalities $x_{2^n+i}=\overline{x_i}$ for each $n\in\NN$ and $i=1,\ld, 2^n$.
That is a peculiar property of self-similarity.

3) The set $x_1,\ld,x_{2^{n+1}}$ can be obtained from the set of binary digits $x_1,\ld,x_{2^n}$
by replacing each digit $0$ with the pair of digits $0,1$ and each $1$ with the pair $1,0$.

4) The continued fraction for the number $\chi$ begins as: 
$$
\chi = [0; 2, 2, 2, 1, 4, 3, 5, 2, 1, 4, 2, 1, 5, 44, 1, ...].
$$

5) The first few convergents of the continued fraction for the number $\chi$ look like this: 
$$
1/2, 2/5, 5/12, 7/17, 33/80, 106/257, 563/1365, 1232/2987, 1795/4352, \ld.
$$
\end{Zam}

\begin{Opr} \label{OprH}
Let's define the mapping $H\colon [0;1)\to[0;1)$ as follows: 
if the point $x\in[0;1)$ has the expansion $x = 0{.}x_1 x_2 \ld$ 
(where periodical $1$ is forbidden for unambiguity of the display), 
then put $H(x) = 0{.}x_1\overline{x_1}x_2\overline{x_2}\ld$.
That is, in order to get $H(x)$, each $0$ is replaced by the pair $01$ and each $1$ by $10$.
\end{Opr}

\begin{Prim} \label{PrimH}
$H(1/2)=7/12$, $H(7/12)=47/80$.
\end{Prim}
{\bf Proof.}
Since the point $1/2$ has the binary expansion $1/2 = 0{.}1000\ld$, then $H(1/2) = 0{.}10\:01\:01\:01\ld = 7/12$.
Furthermore we find: $H(7/12) = 0{.}1001\;0110\;0110\;0110\ld = 47/80$.
\qed

\begin{Predl} \label{PredlSvoistvaH}
1) The function $H$ strictly increases on $[0;1)$.

2) The mapping $H\colon [0;1)\to[0;1)$ has only two invariant points: 
$\chi$ (specified in the definition~\ref{OprCHi}) and $1-\chi$.

3) For any $x\in[0;1/2)$ it is true that $|H^n(x)-\chi| < 1/2^{2^n}$ for each ${n\in\NN}$, 
and $\lim_{n\to\infty}H^n (x)=\chi$.
Similarly, for any ${x\in[1/2;1)}$ it is true that $|H^n(x)-(1-\chi)| < 1/2^{2^n}$ for each $n\in\NN$, 
and $\lim_{n\to\infty}H^n(x)=1-\chi$.

4) The image $H([0;1))$ of the semi-interval $[0;1)$ by mapping $H$ 
lies in the set $E_{1/2}$ of points of the global maximum of the Takagi function $T_{1/2}$ on $[0;1]$.
\end{Predl}
{\bf Proof.}
1) If the number $x=0{.}x_1x_2\ld$ is less than the number $y=0{.}y_1y_2\ld$, 
then the formulae $x_k=y_k$ ($k=1,\ld, n-1$) and $x_n<y_n$ are true for some $n\in\NN$.
Therefore, the number $H(x)=0{.}x_1\overline{x_1}\ld x_n\overline{x_n}x_{n+1}\overline{x_{n+1}}\ld$ is less than the number $H(y)=0{.}y_1\overline{y_1}\ld y_n\overline{y_n}y_{n+1}\overline{y_{n+1}}\ld$.
So the function $H$ strictly increases.

2) Let $x$ be an invariant point of $H$. 
Equating the binary digits of the expansions $x = 0{.}x_1x_2x_3x_4\ld$ 
and $H(x) = 0{.}x_1\overline{x_1}x_2\overline{x_2}\ld$, we get that for any $k\in\NN$, 
the equalities $x_{2k}=\overline{x_k}$ and $x_{2k+1}=x_{k+1}$ are satisfied.
If $x\in[0; 1/2)$, then $x_1=0$, so by virtue of the definition~\ref{OprCHi}, 
the point $x$ coincides with the point $\chi$.
Similarly, if $x\in[1/2;1)$, then $x_1=0$ and $x$ coincides with $1-\chi$.

3) Since any number $x\in[0; 1/2)$ has $x_1=0$, so 
the first $2^n$ binary digits of $H^n(x)$ coincides with the first $2^n$ digits of $\chi$. 
This leads to the inequalities $|H^n(x)-\chi| < 1/2^{2^n}$, 
which lead to the equality $\lim_{n\to\infty}H^n(x)=\chi$.
For $x\in[0; 1/2)$, the reasoning is similar.

4) The statement of item~4) follows from item~2) of the theorem~\ref{TeorPoiskGlobMax}.
\qed

\begin{Teor}
\label{TeorIzmMaksSqrt}
1) Suppose $v\in(1/(2\sqrt2); 1/2)$. 
Then $2v^2\in(1/4; 1/2)$ and the equality $E_v=H(E_{2v^2})$ is fulfilled.
Moreover, the function $H$ (see~definition~\ref{OprH}) bijectively maps the set $E_{2v^2}$ to the set $E_{v}$,
and the inclusion $E_v\subset E_{1/2}$ is true.

2) Suppose $v\in(1/4, 1/2)$. Then $\sqrt{v/2}\in(1/(2\sqrt2); 1/2)$ and the function $H$ bijectively maps the set $E_{v}$ to the set $E_{\sqrt{v/2}}$.
In addition, the inclusion $E_{\sqrt{v/2}}\subset E_{1/2}$ is true.

3) Suppose $v\in(1/2^{1+2^{-n}}; 1/2)$ for some $n\in\{0,1,2,\ld\}$.
Then $q = 2^{2^n-1}v^{2^n}\in(1/4; 1/2)$ and the function $H^n$ bijectively maps the set $E_{q}$ to the set $E_{v}$.
In addition, for any $x\in E_{v}\cap[0;1/2)$ the inequality $|x-\chi|<1/2^{2^n}$ is satisfied.

4) Suppose $v\in(1/4; 1/2)$ and $n\in\{0,1,2,\ld\}$. 
Then $u = 2^{2^{-n}-1}v^{2^{-n}} \in (1/2^{1+2^{-n}}; 1/2)$ and 
the function $H^n$ bijectively maps the set $E_{v}$ to the set $E_{u}$.
In addition, for any $x\in E_{v}\cap[0; 1/2)$ the inequality $|H^n(x)-\chi|<1/2^{2^n}$ is satisfied.

5) Suppose $v\in(1/4; 1/2)$. Then for any point $x\in E_{v}\cap[0;1/2)$ 
the inequalities $|x-\chi| < 2^{1/\log_2(4v^2)} < 2^{-1/(4(1-2v))}$ are satisfied.
\end{Teor}
{\bf Proof.}
1) Let $v\in(1/(2\sqrt2), 1/2)$. Then $2v\in(1/\sqrt2, 1)$.
We consider here only the case when the point $2v$ is consistent with the series $F_{2v}(x) = c_0 + c_1 x + c_2 x^2 + \ldots$.
The case when the point $2v$ is consistent with the point $2v$ is treated similarly.

According to item~1) of the theorem~\ref{TeorSqrtC}, the series $F_{2v}(x)$ has the form $F_{2v}(x) = (1-x) F_{4v^2}(x^2)$, where $F_{4v^2} = b_0 + b_1x + b_2x^2 + \ld$ --- a series consistent with the point $4v^2$.
Equating the coefficients of the series $F_{2v}(x)$ and $(1-x)F_{4v^2}(x^2)$, we get the formulae
\begin{equation} \label{eqc2k2k+1}
c_{2k} = b_k, \quad c_{2k+1} = -b_k \quad \text{when } k=0,1,2,\ld.
\end{equation}
Next, according to item~1a) of the theorem~\ref{TeorPoiskGlobExtr}, the set $E_{v}$ contains two (possibly matching) points $x^-\in[0; 1/2]$ and $x^+=1-x^-\in[1/2;1]$.
The points $x^- $ and $x^+$ can only match if $x^-=x^+=1/2$.
However, due to~\cite [Theorem 4]{Galkin2015}, for $v\in(1/4;1)$, the function $T_v$ does not have a global maximum at $1/2$.
Therefore, the points $x^-$ and $x^+$ are different.
Similarly, since $2v^2\in(1/4, 1/2)$, the set $E_{2v^2}$ also contains two different points $y^-\in[0; 1/2]$ and $y^+\in[1/2;1]$.
It is enough for us to show that $x^- = H(y^-)$ and $x^+ = H(y^+)$.

According to the formula~\eqref{eqx_n^-} we have: 
$x^- = 0{.}x_1x_2\ld$ and $y^- = 0{.}y_1y_2\ld$, where $x_n = (1-c_{n-1})/2$ and $y_n = (1-b_{n-1})/2$ for any $n\in\NN$.
Hence, using the formulae~\eqref{eqc2k2k+1}, we get the equalities 
$x_{2k-1} = (1-c_{2k-2})/2 = (1-b_{k-1})/2 = y_k$ and 
$x_{2k} = (1-c_{2k-1})/2 = (1+b_{k-1})/2 = 1-y_k = \overline{y_k}$ for any $k\in\NN$.
Therefore, by virtue of the definition~\ref{OprH}, we have: 
$x^- = 0{.}x_1x_2x_3x_4\ld = 0{.}y_1\overline{y_1}y_2\overline{y_2}\ld = H(y^-)$.
The equality $x^+ = H(y^+)$ can be proved similarly.

Finally, the inclusion $E_v\subset E_{1/2}$ follows from the proven equality $E_v=H(E_{2v^2})$ 
and item~4) of the proposition~\ref{PredlSvoistvaH}.

2) For the case of $v\in(1/4, 1/2)$, the statement of item~2) follows from the statement of item~1) 
by replacing $v$ with $\sqrt{v/2}$.

3) Let $v\in(1/2^{1+2^{-n}}; 1/2)$, where $n\in\{0,1,2,\ld\}$.
If $n=0$, then the statements being proved are obvious.
If $n\ges1$, then we get that $H^n$ bijectively maps $E_{q}$ to $E_{v}$, by applying statement of item~1) $n$ times.
This implies the equality $E_v = H^n(E_{q})$.
Therefore, by virtue of item~3) of the proposition~\ref{PredlSvoistvaH}, 
the inequality $|x-\chi|<1/2^{2^n}$ is true for any $x\in E_{v}\cap[0;1/2)$.

4) If $v\in(1/4; 1/2)$, then the statement of item~4) of the theorem 
can be derived from the statement of item~3) through replacing $v$ by $2^{2^{-n}-1}v^{2^{-n}}$.

5) If $v\in(1/4; 1/2)$ then set $n = -[\log_2(-\log_2v-1)]$. 
Therefore $n\in\NN$ and $v\in(1/2^{1+2^{-n+1}}; 1/2^{1+2^{-n}}]$.
Hence $1/2^{2^{n-1}} \les 2^{1/\log_2(4v^2)}$.
So, since the inequality $|x-\chi|<1/2^{2^{n-1}}$ is true for any $x\in E_{v}\cap[0;1/2)$ under item~3), 
then we get the required inequalities: $|x - \chi| < 2^{1/\log_2(4v^2)} < 2^{-1/(4(1-2v))}$.
\qed

\medskip
Note that in the case $v = 1/(2\sqrt2)$, the equality $E_v=H(E_{2v^2})$ does not hold.
Indeed, $E_{2v^2} = E_{1/4} = \{1/2\}$ by virtue of item~3) of the theorem~\ref{TeorPoiskGlobMax}, 
so $H(E_{2v^2}) = H(\{1/2\}) = \{5/12\}$ (see example~\ref{PrimH}).
At the same time, $E_v = \{5/12;7/12\}$ (see example~\ref{PrimMaxV1/2sqrt2}).

\begin{Teor}
\label{TeorFormulaChi}
The following formula is correct:
$$
\chi = 1/2 - 1/4\cdot\prod_{n=0}^\infty(1-1/2^{2^n}).
$$
\end{Teor}
{\bf Proof.}
Let's take some number $d\in (1/2;1)$ that is consistent with the series, 
for example $d=2/3$ (see~example~\ref{PrimSoglW2/3}).
Then, by virtue of item~3) of the theorem~\ref{TeorSqrtC}, 
for any $n\in\NN$ the number $w_n = d^{2^{-n}}$ is consistent with the number 
$F_{w_n}(x) = (1-x)(1-x^2)\cdots(1-x^{2^{n-1}}) F_d (x^{2^n})$.
If we put $x=1/2$ here and move to the limit as $n\to\infty$, then we get:
\begin{equation} \label{eqlimFwn}
\lim_{n\to\infty}F_{w_n}(1/2) = \prod_{n=0}^\infty(1-1 / 2^{2^n}) \cdot \lim_{n\to\infty}F_d(1/2^{2^n}) = \prod_{n=0}^\infty(1-1/2^{2^n}).
\end{equation}
According to item~1a) of the theorem~\ref{TeorPoiskGlobExtr}, the set $E_{w_n/2}\cap[0; 1/2]$ consists of a single point $z_n = 1/2 - F_{w_n}(1/2)/4$.
Next, since $w_n/2\in(1/2^{1+2^{-n}}; 1/2)$, then the estimate $|z_n-\chi|<1/2^{2^n}$ is correct, 
according to item~3) of the theorem~\ref{TeorIzmMaksSqrt}.
This implies the equality $\chi =\lim_{n\to\infty}z_n = 1/2 - 1/4\cdot\lim_{n\to\infty}F_{w_n}(1/2)$.
Consequently, taking into account the relation~\eqref{eqlimFwn}, we get the required formula.
\qed 

\subsection{\boldmath Single-ended limits and monotonous decrease of the sets $E_v\cap[0; 1/2]$ by parameter $v$}

Denote by $\widetilde{E}_v$ the set $E_v\cap[0; 1/2]$ consisting of the points of the global maximum of the function $T_v$ on the segment $[0;1/2]$.
\begin{Opr} 
\label{OprLimMa}
Suppose $v_0,y_0\in\RR$.

1) We say that the sets $\widetilde{E}_v$ tend to the point $y_0$ as $v\to v_0-0$, 
and write $\lim_{v\to v_0-0}\widetilde{E}_v = y_0$, 
if $\lim_{v\to v_0-0}\rho(y_0,\widetilde{E}_v) = 0$.

2) The relation $\lim_{v\to v_0+0}\widetilde{E}_v = y_0$ is defined similarly.
\end{Opr}


\begin{Teor}
\label{TeorLimEv}
1) Let $v_0\in[1/4; 1/2)$. Then $\lim_{v\to v_0 - 0}\widetilde{E}_v = \sup\widetilde{E}_{v_0}$ and $\lim_{v\to v_0 + 0}\widetilde{E}_v = \inf\widetilde{E}_{v_0}$. 
In particular, if the set $\widetilde{E}_{v_0}$ consists of a single point 
(when the point $2v_0$ is consistent with a series), 
then $\lim_{v\to v_0 \pm 0}\widetilde{E}_v = \widetilde{E}_{v_0}$.

2) The following estimations are correct:
$1/3=\inf\widetilde{E}_{1/2} <\lim_{v\to 1/2 - 0}\widetilde{E}_v = \chi <\sup \widetilde{E}_{1/2}=5/12$.

3) If $v_0\in[1/4;1/2]$, then the Hausdorff dimension $\dim_H(E_v)$ of the set $E_v$ 
satisfy the equality $\lim_{v\to v_0} \dim_H (E_v) = 0$.
\end{Teor}
{\bf Proof.} 
1) Let ${v_0}\in[1/4; 1/2)$.
Then we have by item~1) of the theorem~\ref{TeorPredelSoglFunc} and by the proposition~\ref{PredlShodUnitFunc}: 
$\lim_{v\to v_0 - 0} F_{2v}^\pm (1/2) = F_{v_0}^-(1/2)$.
In addition, for all $v\in(-1; 1)$ the following equations are true,
due to the formulae \eqref{eqx^-(v)}, \eqref{eqinfEv}, and~\eqref{eqinfEv012} of the theorem~\ref{TeorPoiskGlobExtr}:
\begin{equation} \label{eqInfSupEv}
\inf\widetilde{E}_v = 1/2 - F_{2v}^+(1/2)/4;\qquad
\sup\widetilde{E}_v = 1/2 - F_{2v}^-(1/2)/4.
\end{equation}
Moving $v$ to $v_0-0$, we get: 
$\lim_{v\to v_0 - 0}\inf\widetilde{E}_v = \lim_{v\to v_0 - 0}\sup\widetilde{E}_v = \sup\widetilde{E}_{v_0}$.
Hence the required equality $\lim_{v\to v_0 - 0}\widetilde{E}_v = \sup\widetilde{E}_{v_0}$ is true.

The equality $\lim_{v\to {v_0} + 0}\widetilde{E}_v = \inf\widetilde{E}_u$ can be proved similarly.

2) It follows from item~2) of the theorem~\ref{TeorPredelSoglFunc} and the proposition~\ref{PredlShodUnitFunc},
that $\lim_{w\to 1-0} F_w^{\pm}(1/2) = \prod_{k=0}^{\infty}(1-1/2^{2^k})$. 
Therefore, approaching $v$ to $1/2-0$ in~\eqref{eqInfSupEv}, 
we get: $\lim_{v\to 1/2 - 0}\inf\widetilde{E}_v = \lim_{v\to 1/2 - 0}\sup\widetilde{E}_v = 1/2 - \prod_{k=0}^{\infty}(1-1/2^{2^k})/4$.
Then, this equation and theorem~\ref{TeorFormulaChi} lead us to the required formula $\lim_{v\to 1/2 - 0}\widetilde{E}_v = \chi$.

3) By virtue of the theorem~\ref{TeorPoiskGlobExtr}, for any $v\in(0;1)$ only two cases are possible: 
a) either the point $2v$ is consistent with a series, then $E_v$ contains no more than two points and $\dim_H(E_v) = 0$; 
b) or $2v$ is consistent with a polynomial of degree $N$, then $\dim_H(E_v) = 1/(N+1)$.
Since for any $N\in\NN$ only a finite number of points is consistent with polynomials of degree $N$, 
then $\dim_H (E_v)$ will be arbitrarily small for any $v\neq {v_0}$,
which is close enough to ${v_0}$.
So $\lim_{v\to {v_0}} \dim_H(E_v) = 0$.
\qed

\medskip
The following theorem tells about the monotonicity of sets $\widetilde{E}_v$.
\begin{Teor}
\label{TeorMonotEv}
If $1/4\les u < v < 1/2$, then $\inf\widetilde{E}_u > \sup\widetilde{E}_v$ 
(that is, the set $\widetilde{E}_u$ lies strictly to the right of the set $\widetilde{E}_v$).
\end{Teor}
{\bf Proof.}
By virtue of the theorem~\ref{TeorPoiskGlobExtr}, we have: 
$\inf\widetilde{E}_u = 1/2 - F_{2u}^+(1/2)/4$ and $\sup\widetilde{E}_v = 1/2 - F_{2v}^-(1/2)/4$.
Therefore, we only need to prove the inequality $F_{2u}^+(1/2) < F_{2v}^-(1/2)$.
According to the theorem~\ref{TeorMonotSoglFunc}, the relation $F_{2u}^+ \prec F_{2v}^-$ is true.
Let $c_k^+(u)$ and $c_k (v)^-$ ($k=0,1,2,\ld$) be coefficients of the series $F_{2u}^+$ and $F_{2v}^-$, respectively.
Then, in according with the definition~\ref{OprLeksSravnShod}, 
there exists such $n\in\NN$, that $c_n(u)=-1$, $c_n(v)=1$, and $c_k(u) = c_k (v)$ for all $k=0,1,\ld, n-1$.
So,
\begin{equation} \label{eqF2v-F2u(1/2)}
F_{2v}^-(1/2) - F_{2u}^+(1/2) = \sum_{k=0}^\infty \frac{c_k^-(v)-c_k(u)^+}{2^k} =
\frac{2}{2^n} + \sum_{k=n+1}^\infty \frac{c_k^-(v)-c_k^+(u)}{2^k}.
\end{equation}
Let's show the equality $c_N^-(v)=1$ for some $N>n$ (from which it follows that $c_N^-(v)-c_N^+(u)>-2$).
If the point $2v$ is consistent with the polynomial, 
then this follows from item~1) of the remark~\ref{ZamSoglFun} and item~1) of the definition~\ref{OprSoglFunc}.
If the point $2v$ is consistent with a series, then this follows from the lemma~\ref{LemOdinZnak}, 
because $2v>1/2$ and $1-2v-\ld-(2v)^m < 0$ for some $m\in\NN$.
So then, $c_N^-(v)-c_N^+(u)>-2$ for some $N>n$ and $c_k^-(v)-c_k^+(u)\geqslant-2$ for any $k>n$. 
Whereupon, the desired estimation follows from~\eqref{eqF2v-F2u(1/2)}:
$F_{2v}^+(1/2) - F_{2u}^-(1/2) > 1/2^{n-1} + \sum_{k=n+1}^\infty (-2)/2^k = 0$.
\qed

\subsection{\boldmath Calculation of $v$ by the maximum point of the function $T_v$}

In this section, we present an algorithm of searching for an unknown parameter $v\in (1/4; 1/2)$ 
in the case, when at least one point of the global maximum of the function $T_v$ on the segment $[0;1]$ is known. 
In addition, examples of using this algorithm are provided.

\begin{Teor}[\boldmath algorithm for finding $v$ by the maximum point]
Suppose parameter $v\in (1/4;1/2)$ is unknown, but some point of the global maximum $x_{max}\in[0;1]$ of the function $T_v$ is known.
Then the parameter $v$ can be found using the following algorithm:

Step 1) If $x_{max} > 1/2$, then replace $x_{max}$ with $1-x_{max}\in[0; 1/2]$.

Step 2) Write $x_{max}$ in binary form: $x_{max} = 0{.}x_1x_2x_3\ld$.

Step 3) Build a unitary series $F(x) = c_0 + c_1 x + \ldots + c_n x^n + \ldots$ with coefficients $c_n = 1-2x_{n+1}$, $n=0,1,2,\ld$.
For any $|x|<1$ its sum can also be written as
\begin{equation} \label{eqAlgF}
F(x) = 1/(1-x) - 2(x_1 + x_2\cdot x + \ld + x_{n+1}\cdot x^n + \ld).
\end{equation}

Step 4) Find on the interval $(1/2;1)$ all the roots $x$ of the function $F$,
in which the derivative $F'(x)$ is negative.

Step 5) Select from these roots the only root $w$ for which the series $F$ is intermediate (see definition~\ref{OprSoglFunc}).

Step 6) Calculate the desired parameter by formula $v=w/2$.
\end{Teor}
{\bf Proof.}
Justification of the algorithm follows from the theorems \ref{TeorFw(w)=0}, \ref{TeorProizEdinSogl}, and \ref{TeorPoiskGlobExtr}.
\qed

\begin{Zam}
1) Not every point of the segment $[0; 1/2]$ is the point of the global maximum of the function $T_v$ for some $v\in(-1;1)$. 
For example, it follows from the theorems \ref{TeorPoiskGlobMax}, \ref{TeorIzmMaksSqrt} and theorem~\ref{TeorMonotEv},
that $\bigcup_{v\in(-1; 1)} E_v \subset [1/3, 2/3]$.

2) We can see from the theorems \ref{TeorLimEv} and \ref{TeorMonotEv}, 
that $\bigcup_{v\in(1/4; 1/2)} E_v \subset [\chi, 1-\chi]$.

3) Another inclusion $\bigcup_{v\in[1/(2\sqrt2); 1/2)} \big(E_v\cap[0,1/2]\big) \subset E_{1/2}\cap(\chi,5/12]$
follows from item~1) of the theorem~\ref{TeorIzmMaksSqrt}, item~2) of the theorem~\ref{TeorLimEv}, 
the theorem~\ref{TeorMonotEv} and the example~\ref{PrimMaxV1/2sqrt2}.
This inclusion is not an equality, as examples \ref{PrimXmax317/768} and~\ref{PrimXmax319/768} show.
In the examples \ref{PrimXmax317/768} and~\ref{PrimXmax319/768} we specify points of the set $E_{1/2}\cap(\chi, 5/12]$, 
that are not points of the global maximum of the function $T_v$ on $[0;1]$ for any $v\in[1/(2\sqrt2); 1/2)$.
\end{Zam}

According to the definition~\ref{OprCHi} we have: $\chi = 0{.}0110100110010110\ld$.
In addition, $5/12 = 0{.}01(10)_2$.
Therefore, taking into account Kahane's theorem~\ref{TeorKahane}, we see, that
the set $E_{1/2}\cap(\chi, 5/12]$ includes, in particular, the following four rational points: 
$a_1 = 5/12$, 
$a_2 = 0{.}011010100(101101010010)_2 = 6913/16640 = 0{.}41544471(153846)$, 
$a_3 = 0{.}01101001(10)_2 = 317/768 = 0{.}41276041(6)$ and
$a_4 = 0{.}01101010(01)_2 = 319/768 = 0{.}41536458(3)$.
Let's conjecture, that for any of points $a_i$, $i=1,2,3,4$ there exists some $v_i\in[1/(2\sqrt2); 1/2)$,
for which $a_i$ belongs to the set $E_{v_i}$, $i=1,2,3,4$.
We try to restore four values $v_1$, $v_2$, $v_3$, $v_4$ (if exist) in the following four examples.
\begin{Prim} (cf.~with example~\ref{PrimMaxV1/2sqrt2}).
The point $a_1 = 5/12$ belongs to $E_v$ when $v = 1/(2\sqrt2)$.
\end{Prim}
{\bf Proof.}
Since $a_1 = 0{.}01\;10\;10\;10\ld$, then we have by the formula~\eqref{eqAlgF}: 
$F_1(x) = 1/(1-x) - 2(x + x^2+x^4+x^6+\ld) = (1-2x^2)/(1+x)$.
This function has on $(1/2;1)$ a single root $w=1/\sqrt2$.
The example~\ref{PrimSoglW1/sqrt2} shows that the point $1/\sqrt2$ is consistent with the function $F_1$, 
so $v=w/2=1/(2\sqrt2)$.
\qed

\begin{Prim}
\label{PrimXmax6913/16640}
The point $a_2 = 6913/16640 $ belongs to $E_v$ when $v = w/2 \approx 0{.}381616$, 
where $w\approx0{.}763231$ is the (single) root of the polynomial $1-2x^2+2x^6-2x^8$ on the interval $(1/2;1)$.
\end{Prim}
{\bf Proof.}
Since $a_2 = 0{.}011010100(101101010010)_2$, then we have by the formula~\eqref{eqAlgF}:
$$
F_2(x) = 1/(1-x) - 
2\Bigl(x+x^2+x^4+x^6+(x^9+x^{11}+x^{12}+x^{14}+x^{16}+x^{19})\sum_{k=0}^\infty x^{12k}\Bigr) =
$$
$$ 
= (1-x)(1-x+x^2)(1+x+x^2)(1-2x^2+2x^6-2x^8)/(1+x^6).
$$
On the interval $(1/2;1)$, the polynomial $1-2x^2+2x^6-2x^8 = (1-x^2)^2(1-x^4)-x^8$ 
decreases and takes the values of different signs at the ends, the other cofactors are positive.
This means that the function $F$ has exactly one root $w\approx0{.}763231$ on $(1/2;1)$.
Next, make sure that $w$ is consistent with $F_2$.
\qed

\begin{Prim}
\label{PrimXmax317/768}
The point $a_3 = 317/768$ does not belongs to $E_v$ whatever the parameter $v\in[1/(2\sqrt2); 1/2)$ is.
\end{Prim}
{\bf Proof.}
Since $a_3 = 0{.}01101001(10)_2$, then we have by the formula~\eqref{eqAlgF}:
$$
F_3(x) = 1/(1-x) - 2\Bigl(x+x^2+x^4+x^7 + x^8\sum_{k=0}^\infty x^{2k}\Bigr) =
$$
$$
= 1/(1-x) - 2( x+x^2+x^4+x^7 + x^8/(1-x^2)) = (1-2x^2+2x^6-2x^8)/(1+x).
$$
This function has on $(1/2;1)$ the same single root $w\approx0{.}763231$ 
as the function $F_2$ in the example~\ref{PrimXmax6913/16640}.
But now this root is no consistent with the function $F_3$, since the consistent function $F_2$ is unique.
So, the number $a_3=317/768$ does not belongs to $E_v$ whatever the parameter $v\in[1/(2\sqrt2); 1/2)$ is.
\qed

\begin{Prim}
\label{PrimXmax319/768}
The point $a_4 = 319/768$ does not belongs to $E_v$ whatever the parameter $v\in[1/(2\sqrt2); 1/2)$ is.
\end{Prim}
{\bf Proof.}
Since $a_4 = 0{.}01101010(01)_2$, then we have by the formula~\eqref{eqAlgF}:
$$
F_4(x) = 1/(1-x) - 2\Bigl(x+x^2+x^4+x^6 + x^9\sum_{k=0}^\infty x^{2k}\Bigl) = (1-2x^2+2x^8)/(1+x).
$$
Where as $(1-2x^2+2x^8) = ((4+4x^2+3x^4)(2-3x^2)^2+5x^8)/16 > 0$ for any $x\in\RR$,
then $F_4$ has no real roots.
So, the number $a_4=319/768$ also does not belongs to $E_v$ whatever the parameter $v\in[1/(2\sqrt2); 1/2)$ is.
\qed

\subsection{Several examples of calculating of points of global maximum for functions in exponential Takagi class}
\label{PrimerMaks}

\begin{Prim} 
For any $v$ in the interval $-1 < v < -(1+\sqrt5)/4\approx -0{.}809017$, 
the global maximum of the function $T_v$ on the segment $[0;1]$ equals $(2+v)/(5(1-v^2))$ 
and is reached at two points: $2/5$ and $3/5$.
If $v = -(1+\sqrt5)/4$, then this maximum equals $(15+\sqrt5)/25$, 
and is reached at four points: $2/5$, $19/40$, $21/40$ and $3/5$.
\end{Prim}
{\bf Proof.}
Let's use item~4) of the theorem~\ref{TeorPoiskGlobMax} for $k=1$. If follows from it, that:

1) the polynomial $\widetilde{P}_2 (x) = 1-2x-4x^2$ has a single negative root $v_1$
(it is clear that $v_1 = -(1+\sqrt5)/4$);

2) in the case $-1 < v < -(1+\sqrt5)/4$, the function $T_v$ reaches a global maximum at two points on $[0;1]$: 
$x_{max,1}^- = 2/5$ and $x_{max,1}^+ = 3/5$;

3) in the case $v = -(1+\sqrt5)/4$, four maximum points exists:
$x_{max,1}^\pm$, and in addition $x_{max,2}^- = 19/40$, $x_{max, 2}^+ = 21/40$.

The global maximum value can be obtained from the formula~\eqref{eqmaxTvk}.
\qed

\begin{Prim}[maximum points for the sequence $\{T_{v_n}\}$]
For any $n\ges 2$ the polynomial $\widetilde{P}_n (x) = 1-2x-4x^2-\ld-2^nx^n$ has a single positive root $v_n$.
The sequence $\{v_n\}$ lies in the interval $(1/4;1/(2\sqrt{2}))$ and is strictly decreasing. 
Moreover, the following asymptotic formula is correct:
\begin{equation} \label{eqAsimptVn}
v_n = 1/4+1/2^{n+3}+\underline{O}(n / 2^{2n}) \quad ({n\to\infty}).
\end{equation}
For any $n\ges 2$, the set $E_{v_n}$ has a Hausdorff dimension $1/(n+1)$ and consists of all points $x_{max}$ with a binary expansion of the form
\begin{equation} \label{eqTo4kMaksVn}
x_{max} = 0{.}
\left[\begin{array}{l}011\ld1\\\text{or}\\100\ld0\\\end{array}\right.
\left[\begin{array}{l}011\ld1\\\text{or}\\100\ld0\\\end{array}\right.
\ld,
\end{equation}
where each block of ones and zeros contains $n$ digits.
In addition, two more asymptotic equalities hold:
\begin{equation} \label{eqinfEVn}
\inf E_{v_n} = 1/2 - 1/(2^{n+2}-2),
\end{equation}
\begin{equation} \label{eqMVn}
M_{v_n} = 1/(2(1-v_n^{n+1})) = 1/2+1/2^{2n+3}+\underline{O}(n/2^{3n}) \quad ({n\to\infty}).
\end{equation}
\end{Prim}
{\bf Proof.}
The uniqueness of the positive root $v_n$ and the asymptotic formula~\eqref{eqAsimptVn} 
follow from items 1) and~2) of the theorem~\ref{TeorEdinPolozKor}, with $v_n = u_n/2$.

Note, that for every $n\ges 2$ the consistency of the polynomial $P_n(x) = 1-x-x^2-\ld-x^n$ with the number $2v_n$ 
follows from item~3) of the theorem~\ref{TeorEdinPolozKor}. 

Next, use item~2) of the theorem~\ref{TeorPoiskGlobExtr}.
The equality $P_{n}(1/2) = 1/2^n$ implies 
the equalities $a_n (v_n) = 1/2 - P_{n}(1/2)/4 - 1/2^{n+2} = 0{.}011\ld1_2$ and $b_n (v_n) = 1/2 - P_{n}(1/2)/4 + 1/2^{n+2} = 1/2$, 
as well as the formulae \eqref{eqTo4kMaksVn} (for maximum points) 
and \eqref{eqinfEVn} for ($\inf E_v$).

Due to the equality $(1-2v_n)\widetilde{P}_n(v_n) = 1-4v_n+2^{n+1}v_n^{n+1} = 0$ and the formula~\eqref{eqAsimptVn}, 
we get: $v_n^{n+1} = (4v_n-1)/2^{n+1} = 1/2^{2n+2}+\underline{O}(n/2^{3n}) \quad ({n\to\infty})$.
Hence, by virtue of the equality $M_{v_n, n} = S_{v, n}(b_n(v_n)) = S_{v, n}(1/2) = 1/2$ and the formula~\eqref{eqMv_MvN}, 
the required relation~\eqref{eqMVn} follows.
\qed

\begin{Prim} 
If $v = 1/(2\sqrt[4]{2}) \approx 0{.}420448$, then the function $T_v$ has two maximum points on the segment $[0,1]$: 
$33/80 = 0{.}4125$ and $47/80 = 0{.}5875$.
\end{Prim}
{\bf Proof.}
In the example~\ref{PrimSoglW_sqrt4_2}, it is shown that the point $2v = 1/\sqrt[4]{2}$ is consistent with the series $F_{2v}$,
which has the sum $F_{2v}(x) = (1-x)(1-2x^4)/(1+x^2)$ at $|x|<1$.
Then, using the formulae \eqref{eqx^-(v)} and~\eqref{eqx^+(v)}, we get the maximum points of the function $T_v$ on $[0, 1]$: 
$x_{max}^- = 1/2-1/4\cdot F_{2v}(1/2) = 33/80$ and $x_{max}^+ = 1-x_{max}^- = 47/80$.
\qed

\begin{Prim} 
If $v = 1/3 \approx 0{.}3333$, then the function $T_v$ has two maximum points on the segment $[0,1]$: 
$x_{max}^- \approx 0{.}419240$ and $x_{max}^+ \approx 0{.}580760$.
\end{Prim}
{\bf Proof.}
In the example~\ref{PrimSoglW2/3}, it is shown that the point $2v = 2/3$ is consistent with the series $F$,
and the first few terms of $F$ is given.
Knowing these members and using the formulae \eqref{eqx^-(v)} with~\eqref{eqx^+(v)}, 
it is possible to approximate the maximum points of the function $T_v$ on $[0,1]$: 
$x_{max}^- = 1/2 - 1/4\cdot F(1/2) \approx 0{.}419240$ and $x_{max}^+ = 1-x_{max}^- \approx 0{.}580760$.
\qed

\begin{Prim} 
If $v = 1/(2\sqrt{3})\approx 0{.}288675$, then the function $T_v$ has two maximum points on the segment $[0,1]$: 
$x_{max}^- \approx 0{.}455393$ and $x_{max}^+ \approx 0{.}544607$.
\end{Prim}
{\bf Proof}
This fact follows, similar to the previous example, from the formulae \eqref{eqx^-(v)} with~\eqref{eqx^+(v)}, 
and from the example~\ref{PrimSoglW1/sqrt3}, 
that shows the initial terms of the consistent series for the point $2v = 1/\sqrt{3}$.
\qed





\begin{thebibliography}{99}

\bibitem{HanSchied2020} 
Han~X., Schied~A. {\it Step roots of Littlewood polynomials and the extrema of functions in the Takagi class}. //
arXiv:2001.01348v2 [math.CA]

\bibitem{MishuraSchied2019} 
Mishura~Y., Schied~Y. {\it On (signed) Takagi-Landsberg functions: pth variation, maximum, and modulus of continuity}. //
J. Math. Anal. Appl. Vol.~473. 2019. No.~1. Pp.~258--272.

\bibitem{Han2019} 
Han~X. {\it On the extrema of functions in the Takagi class}. //
Master's thesis. 2019. University of Waterloo, Ontario, Canada. 

\bibitem{Galkin2015}
Galkin~O.~E., Galkina~S.~Yu. {\it On properties of functions in exponential Takagi class}. //
Ufa Mathematical Journal. Vol.~7. 2015. No.~3. Pp.~29--38.

\bibitem{Mandelbrot}
B.~B.~Mandelbrot, {\it Fractal landscapes without creases and with rivers}. // 
Appendix A in {\it The Science of Fractal Images} H.-O.~Peitgen, D.~Saupe Ed. Springer-Verlag. N.~Y. 1988.

\bibitem{Takagi}
Takagi~T. {\it A simple example of the continuous function without derivative}. //
Tokyo Sugaku-Butsurigakkwai Hokoku. Vol.~1. 1901. Pp. 176--177.
https://doi.org/10.11429/subutsuhokoku1901.1.F176


\bibitem{HataYamaguti84} 
Hata~M., Yamaguti~M. {\it Takagi function and its generalization}. //
Japan J. Appl. Math. Vol.~1. 1984. Pp.~183–-199.

\bibitem{HataYamaguti83}
Hata~M., Yamaguti~M. {\it Weierstrass’s function and chaos}. //
Hokkaido Math. J. Vol.~12. 1983. Pp.~333–-342.

\bibitem{Kono} 
K\^ono~N. {\it On generalized Takagi functions}. //
Acta Math. Hungar. Vol.~49. 1987. Pp.~315–-324.

\bibitem{Baba} 
Baba~Y. {\it On maxima of Takagi-van derWaerden functions}. //
Proc. Amer. Math. Soc. Vol.~91. 1984. No.~3. Pp.~373–-376.

\bibitem{Kahane} 
Kahane~J.-P. {\it Sur l’exemple, donn\'e par M. de Rham, d’une fonction continue sans d\'eriv\'ee}. //
Enseignement Math. Vol.~5. 1959. Pp.~53-–57.

\bibitem{AllaartKawamuraSurv} 
Allaart~P.~C., Kawamura~K. {\it The Takagi function: a survey}. //
Real Anal. Exchange. Vol.~37. 2011/12. No.~1. Pp.~1--54.

\bibitem{AllaartKawamura} 
Allaart~P.~C., Kawamura~K. 
{\it Extreme values of some continuous, nowhere differentiable functions}. //
Math. Proc. Camb. Phil. Soc. Vol.~140. 2006. No.~2. Pp.~269-–295.

\bibitem{Lagarias} 
Lagarias~J.~C. {\it The Takagi Function and Its Properties} //
In: Functions and Number Theory and Their Probabilistic Aspects, 
RIMS Kokyuroku Bessatsu B34. Aug. 2012. Pp.~153--189.

\bibitem{Martynov}
Martynov~B. {\it On maxima of the van der Waerden function}. // 
Kvant. 1982. No.~6. Pp.~8-14.

\bibitem{Galkina}
Galkina~S.~Yu. {\it Fourier-Haar coefficients of functions of bounded variation}. // 
Mat. Zametki. Vol.~51. 1992. No.~1. Pp.~42-54.
	
\bibitem{Warden}
Warden~В.~L. van der {\it Ein einfaches Beispiel einer nicht differenzirenbaren stetigen Function}. // 
Math. Zeitschrift. 1930. No.~32. Pp.~474--475.

\bibitem{Ledrappier}
Ledrappier~F. {\it On the dimension of some graphs}. //
Contemp. Math. Vol.~135. 1992. Pp.~285–-293.

\bibitem{Tabor1}
Tabor~J., Tabor~J. {\it Generalized approximate midconvexity}. //
Control Cybernet. Vol.~38. 2009. No.~3. Pp.~655-–669.

\bibitem{Tabor2} 
Tabor~J., Tabor~J. {\it Takagi functions and approximate midconvexity}. //
J. Math. Anal. Appl. Vol.~356. 2009. No.~2. Pp.~729–-737.

\bibitem{LagariasMaddock}
Lagarias~J.~C., Maddock~Z. {\it Level sets of the Takagi function: local level sets}. //
Monat. Math. Vol.~166. 2012. No.~2. Pp.~201--238.

\bibitem{Solomyak}
Solomyak~B. {\it On the random series $\sum\pm\lambda^n$ (an Erd\"os problem)}. //
Ann. Math. Vol.~142. 1995. No.~3. Pp.~611-–625.

\bibitem{Privalov}
Privalov~I.~I. {\it Introduction to the theory of functions of a complex variable}
Moscow. Nauka. 1984. 432~p.

\bibitem{Spurrier}
Spurrier~K.~G. {\it Continuous Nowhere Differentiable Functions} 
Senior Thesis. South Carolina Honors College. April 2004.

\bibitem{ShidfarSabetfakhri}
Shidfar~A., Sabetfakhri K. 
{\it On the Continuity of Van Der Waerden’s Function in the	H\"older Sense}. //
Amer. Math. Monthly. Vol.~93. 1986. No.~5. Pp.~375–-376.

\bibitem{deRham}
de Rham G. {\it Sur un exemple de fonction continue sans derivee} 
Enseign. Math. Vol.~3. 1957. Pp.~71-–72.

\bibitem{Allaart}
Allaart~P.~C. {\it How large are the level sets of the Takagi function?}. //
Monatshefte f\"ur Mathematik. Vol.~167. 2012. No.~3-4. Pp.~311--331.

\bibitem{Kostrikin}
Kostrikin~A.~I. {\it Introduction to algebra. Part I. Basics of algebra}. 
M. Fizmatlit. 2004. 272~p.

\bibitem{Kruppel2007} 
Kr\"uppel~M. {\it On the extrema and the improper derivatives of Takagi’s continuous nowhere differentiable function}. //
Rostock. Math. Kolloq. Vol.~62. 2007. Pp.~41-–59.

\bibitem{Kruppel2008}
Kr\"uppel~M. {\it Takagi’s continuous nowhere differentiable function and binary digital sums}. //
Rostock. Math. Kolloq. Vol.~63. 2008. Pp.~37--54.

\bibitem{Trollope}
Trollope~J.~R. {\it An explicit expression for binary digital sums}. //
Math. Mag. Vol.~41. 1968. Pp.~21--25.

\end{thebibliography}
\end{document}